\documentclass[12pt]{amsart}

\usepackage[T1]{fontenc}
\usepackage[hmargin=1 in, vmargin = 1 in]{geometry}
\usepackage{libertinus}

\usepackage{microtype}
\usepackage{amssymb, amsmath, amsfonts, amsthm}
\usepackage{enumitem}
\usepackage{aliascnt}

\usepackage{hyperref}

\hypersetup{
  hidelinks,
  pdftitle={Invariant Algebraic Connections on Connected Reductive Groups},
  pdfauthor={Rudrendra Kashyap and Ruoxi Li},
  pdfsubject={Algebraic connections on connected reductive groups},
  pdfkeywords={algebraic connections, regular-singular connections, reductive groups,
  finite central descent, idempotent completion, algebraic D-modules}
}

\numberwithin{equation}{section}
\setlist[enumerate]{label=(\arabic*),leftmargin=2.2em}

\newtheorem{theorem}{Theorem}[section]

\newaliascnt{proposition}{theorem}
\newtheorem{proposition}[proposition]{Proposition}
\aliascntresetthe{proposition}

\newaliascnt{corollary}{theorem}
\newtheorem{corollary}[corollary]{Corollary}
\aliascntresetthe{corollary}

\newaliascnt{lemma}{theorem}
\newtheorem{lemma}[lemma]{Lemma}
\aliascntresetthe{lemma}

\newaliascnt{conjecture}{theorem}
\newtheorem{conjecture}[conjecture]{Conjecture}
\aliascntresetthe{conjecture}

\theoremstyle{definition}
\newaliascnt{definition}{theorem}
\newtheorem{definition}[definition]{Definition}
\aliascntresetthe{definition}

\newaliascnt{example}{theorem}
\newtheorem{example}[example]{Example}
\aliascntresetthe{example}

\newaliascnt{convention}{theorem}
\newtheorem{convention}[convention]{Convention}
\aliascntresetthe{convention}

\theoremstyle{remark}
\newaliascnt{remark}{theorem}
\newtheorem{remark}[remark]{Remark}
\aliascntresetthe{remark}

\usepackage[capitalize,noabbrev]{cleveref}
\crefname{theorem}{Theorem}{Theorems}
\Crefname{theorem}{Theorem}{Theorems}
\crefname{proposition}{Proposition}{Propositions}
\Crefname{proposition}{Proposition}{Propositions}
\crefname{corollary}{Corollary}{Corollaries}
\Crefname{corollary}{Corollary}{Corollaries}
\crefname{lemma}{Lemma}{Lemmas}
\Crefname{lemma}{Lemma}{Lemmas}
\crefname{conjecture}{Conjecture}{Conjectures}
\Crefname{conjecture}{Conjecture}{Conjectures}
\crefname{definition}{Definition}{Definitions}
\Crefname{definition}{Definition}{Definitions}
\crefname{example}{Example}{Examples}
\Crefname{example}{Example}{Examples}
\crefname{remark}{Remark}{Remarks}
\Crefname{remark}{Remark}{Remarks}
\crefname{convention}{Convention}{Conventions}
\Crefname{convention}{Convention}{Conventions}

\AddToHook{env/theorem/begin}{\crefalias{section}{theorem}}
\AddToHook{env/proposition/begin}{\crefalias{section}{proposition}}
\AddToHook{env/corollary/begin}{\crefalias{section}{corollary}}
\AddToHook{env/lemma/begin}{\crefalias{section}{lemma}}
\AddToHook{env/conjecture/begin}{\crefalias{section}{conjecture}}
\AddToHook{env/definition/begin}{\crefalias{section}{definition}}
\AddToHook{env/example/begin}{\crefalias{section}{example}}
\AddToHook{env/convention/begin}{\crefalias{section}{convention}}
\AddToHook{env/remark/begin}{\crefalias{section}{remark}}

\DeclareMathOperator{\Det}{det}
\DeclareMathOperator{\GL}{GL}
\DeclareMathOperator{\SL}{SL}
\DeclareMathOperator{\PGL}{PGL}
\DeclareMathOperator{\Spec}{Spec}
\DeclareMathOperator{\Lie}{Lie}
\DeclareMathOperator{\Hom}{Hom}
\DeclareMathOperator{\End}{End}
\DeclareMathOperator{\Aut}{Aut}
\DeclareMathOperator{\Ext}{Ext}

\DeclareMathOperator{\pr}{pr}

\DeclareMathOperator{\diag}{diag}
\DeclareMathOperator{\Kar}{Kar}

\DeclareMathOperator{\Der}{Der}
\DeclareMathOperator{\DR}{DR}

\newcommand{\C}{\mathbb C}
\newcommand{\Gm}{\mathbb G_m}
\newcommand{\cO}{\mathcal O}
\newcommand{\g}{\mathfrak g}

\newcommand{\Invtriv}{\operatorname{Conn}^{\mathrm{triv}}}
\newcommand{\Invdesc}{\operatorname{Conn}^{\mathrm{desc}}}
\newcommand{\Connrs}{\operatorname{Conn}_{\mathrm{rs}}}
\newcommand{\LocSysC}{\operatorname{LocSys}}
\newcommand{\RepC}{\operatorname{Rep}}
\newcommand{\Gab}{G^{\mathrm{ab}}}
\newcommand{\Gder}{G^{\mathrm{der}}}
\newcommand{\Gsc}{G^{\mathrm{sc}}}
\newcommand{\HdR}{H_{\mathrm{dR}}}
\newcommand{\RGammadR}{R\Gamma_{\mathrm{dR}}}
\newcommand{\xto}[1]{\xrightarrow{#1}}
\newcommand{\isoto}{\xrightarrow{\sim}}

\title[Invariant Algebraic Connections on Reductive Groups]{Invariant Algebraic Connections on Connected Reductive Groups}
\author{Rudrendra Kashyap}
\address{Department of Mathematics, University of Pittsburgh, Pittsburgh, PA, USA}
\email{ruk26@pitt.edu}
\author{Ruoxi Li}
\address{Department of Mathematics, University of California, Riverside, Riverside, CA, USA}
\email{ruoxi.li1@ucr.edu}


\begin{document}

\begin{abstract}
Let $G$ be a connected complex reductive algebraic group. We study finite-rank flat algebraic connections on trivial vector bundles whose connection forms are left invariant, allowing arbitrary algebraic horizontal morphisms.

We prove that a flat algebraic connection on $G$ is regular-singular
if and only if its pullback to the canonical finite central cover is
isomorphic to a left-invariant flat algebraic connection on a trivial vector
bundle. Moreover, every regular-singular algebraic connection on $G$
is a direct summand of such a connection.

We characterize the essential image of pullback from the
abelianization by the vanishing of a derived-monodromy obstruction and show that pullback is an equivalence precisely when
$\Gder$ is simply connected.
For semisimple
$G$, regular-singular connections are classified by
finite-dimensional representations of the finite central kernel of
the simply connected cover. We also classify regular-singular
connections on $\GL_r$ by pullback along the determinant, give a
$\PGL_2$ counterexample to the abelianization classification of left-invariant trivial-bundle algebraic connections, and compute de Rham and
Betti cohomology in the semisimple case.
\end{abstract}

\maketitle
\setcounter{tocdepth}{1}

\tableofcontents
\section{Introduction}

Throughout this paper, we work over the field of complex numbers $\C$, and $G$ denotes a connected reductive algebraic group over $\C$. A left-invariant flat algebraic connection on a trivial vector bundle is encoded infinitesimally by a representation of $\Lie(G)$. On the simply connected semisimple factor of the canonical cover, that representation integrates and can be removed by an algebraic gauge. What remains is carried by the central torus together with finite descent data. This separation between infinitesimal data and global topology is the organizing principle of the paper.

This formulation is related to the classical theory of invariant connections on homogeneous bundles initiated by Wang~\cite{Wan58}, but the category considered here is different: we impose flatness in the algebraic category and allow arbitrary algebraic horizontal morphisms. Wei~\cite{Wei25} formulated the corresponding free-module problem, established several basic cases, and proposed the abelianization conjecture recalled in \cref{conj:wei}. 

The principal question is to determine which algebraic connections are generated by left-invariant algebraic connections on trivial vector bundles after finite central descent and passage to direct summands. This paper gives an intrinsic description of the regular-singular de Rham category of $G$ and identifies the precise obstruction to pullback from the abelianization.

The ambient regular-singular category has a classical description: Deligne's regular Riemann--Hilbert correspondence identifies regular-singular algebraic connections with complex local systems~\cite[II.5.9]{Del70}.
This paper determines which monodromy representations are realized by left-invariant algebraic connection forms on
trivial vector bundles, how finite central descent enlarges this class, and why passage to direct summands recovers the entire regular-singular category.

We study left-invariant flat algebraic connections on trivial vector bundles, with equivalence given by arbitrary algebraic horizontal isomorphisms. For a connected reductive group, finite central descent can produce a nontrivial vector bundle downstairs. Writing $\Gder=[G,G]$, this leads to a second category, $\Invdesc(G)$, defined using the canonical central isogeny
$
  (\Gder)^{\mathrm{sc}}\times Z(G)^0\longrightarrow G.
$

The canonical cover separates the two sources of information in the problem. The semisimple part is controlled by integration of Lie algebra representations and algebraic gauge transformations; the central torus carries the continuous monodromy; and the finite kernel records the residual global descent obstruction. The resulting classification is therefore governed by the fundamental group and the root datum.

There is also a link with the invariant-theoretic problem for commuting schemes. Let $H$ be a reductive algebraic group with Lie algebra $\mathfrak h$, let $\mathfrak t\subset\mathfrak h$ be the Lie
algebra of a maximal torus, and let $W$ be the corresponding Weyl group. For $d\geq1$, the $d$-fold commuting scheme is the closed subscheme
\[
  \mathcal C_{\mathfrak h}^{d}
  :=
  \left\{
    (x_1,\ldots,x_d)\in\mathfrak h^{d}
    \ \middle|\
    [x_i,x_j]=0\text{ for all }i,j
  \right\}.
\]
Restriction along
$\mathfrak t^{d}\hookrightarrow\mathcal C_{\mathfrak h}^{d}$ induces a natural homomorphism
\[
  \vartheta_{H,d}:
  \C[\mathcal C_{\mathfrak h}^{d}]^{H}
  \longrightarrow
  \C[\mathfrak t^{d}]^{W}.
\]
The higher-dimensional Chevalley restriction problem asks whether $\vartheta_{H,d}$ is an isomorphism. For $d=1$, this is the classical Chevalley restriction theorem. The higher-dimensional statement is
known for general linear, symplectic, and orthogonal
groups~\cite{Vac07,CN21,SXX23}, but the full scheme-theoretic problem remains open for general reductive groups. For $H=\GL_n$, the same commuting scheme also appears in the present setting: if $T\cong(\Gm)^d$, then, by \cref{prop:invariant-forms-lie-reps}, a left-invariant flat connection on the trivial rank-$n$ bundle over $T$ has the form
\[
  d+\sum_{i=1}^{d}A_i\frac{dz_i}{z_i},
  \qquad [A_i,A_j]=0,
\]
and is therefore determined by a point of
$\mathcal C_{\mathfrak{gl}_n}^{d}$. The invariant-theoretic problem studies the ring of functions on this scheme that are invariant under simultaneous conjugation, while the present paper considers the
corresponding connections up to arbitrary algebraic horizontal isomorphism.

\begin{convention}
For a category $\mathcal{C}$, we denote by $[\mathcal{C}]$ the set of isomorphism classes of its objects. We use $=$ for canonical
identifications of sets and for full subcategories of a fixed ambient category that coincide, $\cong$ for isomorphisms of objects, and $\simeq$ for equivalences of categories.    
\end{convention}
\subsection*{Organization}
\Cref{sec:main-results} collects the precise statements of main results.
\Cref{sec:prelim} defines the two categories, records the gauge convention, gives a dictionary with Wei's notation, proves the torus classification directly, and proves the regular-singular characterization. \Cref{sec:semisimple} treats semisimple groups and the $\PGL_2$ counterexample. \Cref{sec:glr} gives both the topological and explicit algebraic determinant arguments. \Cref{sec:reductive} proves the root-datum classification, the derived-monodromy theorem, and the general Karoubi-completion theorem. \Cref{sec:derham,sec:local-systems} establish the de Rham and Betti consequences in the semisimple case.
\subsection*{Acknowledgements}
Both authors were partially supported by NSF grant DMS-2402553.
\section{Main results}\label{sec:main-results}
\subsection{The two categories and regular singularity}

We begin with the definitions. An algebraic connection belongs to $\Invtriv(G)$ if its underlying
vector bundle is algebraically trivial and, in some algebraic
trivialization, its connection form is left invariant. We take
$\Invtriv(G)$ to be full and replete. It belongs to $\Invdesc(G)$ if
its pullback along the canonical central cover
\[
  (\Gder)^{\mathrm{sc}}\times Z(G)^0\longrightarrow G
\]
belongs to
\[
  \Invtriv\bigl((\Gder)^{\mathrm{sc}}\times Z(G)^0\bigr),
\]
where $\Gder=[G,G]$ and $Z(G)^0$ is the connected center.

For $n\geq1$, we write $\Invtriv_n(G)$ and $\Invdesc_n(G)$ for the full subcategories of rank-$n$ objects. For an additive category $\mathcal A$, $\Kar(\mathcal A)$ denotes its idempotent completion, also called its Karoubi envelope.

Our basic structural theorem identifies the second category intrinsically.

\begin{theorem}[\cref{thm:rs}]\label{thm:intro-rs}
For every connected complex reductive group $G$,
\[
  \Invdesc(G)=\Connrs(G)
\]
as full subcategories of the category of finite-rank algebraic connections on $G$. Horizontal sections and monodromy therefore give exact tensor equivalences
\[
  \Invdesc(G)\isoto \LocSysC(G^{\mathrm{an}})
  \isoto \RepC\bigl(\pi_1(G^{\mathrm{an}},e)\bigr).
\]
\end{theorem}
If $T\subset G$ is a maximal torus and $Q^\vee\subset X_*(T)$ is the coroot lattice, the standard identification of the topological and algebraic fundamental groups gives the following concrete description.
\begin{corollary}[\cref{cor:root-classification}]
For every maximal torus $T\subset G$, there is an exact tensor equivalence
\[
  \Invdesc(G)\isoto\RepC\bigl(X_*(T)/Q^\vee\bigr).
\]
\end{corollary}
The fundamental-group exact sequence is established in \cref{prop:pi-exact}, and the root-datum description and categorical classification are proved in \cref{prop:root-datum,cor:root-classification}.

\subsection{The semisimple case}

Let $G$ be a connected semisimple group, let $\pi:\Gsc\twoheadrightarrow G$ be the simply connected central cover, and put $\Gamma=\ker(\pi)$. The two categories have different classifications.

\begin{proposition}[\cref{prop:semisimple-triv}]\label{prop:intro-sem-simple-triv}
For every $n\geq1$, restriction to $\Gamma$ induces a natural bijection
\[
  [\Invtriv_n(G)]
  =
  \operatorname{im}\!\left(
  \frac{\Hom_{\mathrm{alg.grp}}(\Gsc,\GL_n)}{\GL_n}
  \longrightarrow
  \frac{\Hom(\Gamma,\GL_n)}{\GL_n}
  \right).
\]
\end{proposition}

\begin{theorem}[\cref{thm:semisimple-desc} and \cref{cor:semisimple-kar}]\label{thm:intro-sem-desc}
There is an exact tensor equivalence
\[
  \Invdesc(G)\isoto \RepC(\Gamma).
\]
In particular,
\[
  [\Invdesc_n(G)]
  =\Hom(\Gamma,\GL_n)/\GL_n.
\]
\end{theorem}

For $G=\PGL_2$, the standard representation of $\SL_2$ produces a nontrivial rank-two object of $\Invtriv(\PGL_2)$, contradicting \cref{conj:wei}, because $\PGL_2^{\mathrm{ab}}=\{1\}$. The nontrivial character of $\mu_2$ produces, by contrast, a rank-one object only after finite descent. The construction and the descended sign line are given in \cref{ex:std,ex:sign-line}.

\subsection{The general linear group}

Let $\Det:\GL_r\to\Gm$ be the determinant.

\begin{theorem}[\cref{thm:det-equivalence}]\label{thm:intro-gl}
Determinant pullback induces exact tensor equivalences
\[
  \Det^*: \Invdesc(\Gm)\isoto\Invdesc(\GL_r),
  \qquad
  \Det^*: \Invtriv(\Gm)\isoto\Invtriv(\GL_r).
\]
In particular,
\[
  [\Invdesc_n(\GL_r)]
 =\Hom(\mathbb Z,\GL_n)/\GL_n.
\]
\end{theorem}
A homomorphism $\mathbb Z\to\GL_n$ is determined by the image of
$1$. Hence the displayed quotient is the set of conjugacy classes of
invertible $n\times n$ matrices, equivalently their Jordan normal forms.
\cref{sec:glr} also contains the explicit algebraic gauge calculation.

\subsection{Reductive groups and Karoubi completion}

Put $\Gamma=\pi_1((\Gder)^{\mathrm{an}},e)$, and let $\mathrm{ab}:G\to\Gab$ be the abelianization. The topology of $G$ is organized by the short exact sequence
\[
  0\longrightarrow\Gamma\longrightarrow\pi_1(G^{\mathrm{an}},e)
  \longrightarrow\pi_1((\Gab)^{\mathrm{an}},e)\longrightarrow0,
\]
proved in \cref{prop:pi-exact}.
Restriction of monodromy gives an exact tensor functor
\[
  \operatorname{Res}_{\mathrm{der}}:\Invdesc(G)\longrightarrow \RepC(\Gamma).
\]
\begin{theorem}[\cref{thm:abelianization}]\label{thm:intro-ab}
Pullback along $\mathrm{ab}:G\to\Gab$ induces a fully faithful exact
tensor functor
\[
  \mathrm{ab}^*:
  \Invdesc(\Gab)\longrightarrow\Invdesc(G).
\]
Its essential image is the full subcategory consisting precisely of those objects whose
derived monodromy is trivial:
\[
  \left\{
    M\in\Invdesc(G)
    \ \middle|\
    \operatorname{Res}_{\mathrm{der}}(M)
    \cong\mathbf 1_\Gamma^{\oplus\operatorname{rank}(M)}
  \right\}.
\]
Consequently, for every $n\geq1$, pullback induces an injection
\[
  [\Invdesc_n(\Gab)]
  \longrightarrow
  [\Invdesc_n(G)].
\]
These maps are surjective for every $n\geq1$ if and only if
$\Gamma=1$, equivalently, if and only if $\Gder$ is simply
connected.
\end{theorem}

An algebraic construction of the derived-monodromy functor is given in \cref{subsec:algebraic-derived-monodromy}.

The relation between the two categories is uniform for all connected reductive groups.

\begin{theorem}[\cref{thm:general-kar}]\label{thm:intro-kar}
The natural inclusion induces a tensor equivalence
\[
  \Kar\bigl(\Invtriv(G)\bigr)\isoto \Invdesc(G).
\]
Equivalently, every regular-singular algebraic connection on $G$ is a direct summand of a left-invariant algebraic connection on a trivial vector bundle.
\end{theorem}
\begin{corollary}[\cref{cor:simples-semisimple}]
Every simple object of $\Invdesc(G)$ has rank one. Moreover, this category is semisimple if and only if $G$ is semisimple.
\end{corollary}
\subsection{The Wei--Regular-Singular Criterion}
The following result gives a group-theoretic criterion for Wei's conjecture (see \cref{conj:wei}) to hold in every rank and, equivalently, for invariant trivial-bundle connections to exhaust all regular-singular connections.
\begin{theorem}[\cref{cor:wei-iff}]
For a connected complex reductive group $G$, the following conditions are equivalent:
\begin{enumerate}
  \item $\Gder$ is simply connected;
  \item pullback along $\mathrm{ab}:G\to\Gab$ induces an equivalence
  \[
    \mathrm{ab}^*:\Invtriv(\Gab)\isoto\Invtriv(G);
  \]
  \item Wei's \cref{conj:wei} holds for $G$ in every rank;
  \item $\Invtriv(G)$ is idempotent complete;
  \item $\Invtriv(G)=\Invdesc(G)=\Connrs(G)$.
\end{enumerate}
\end{theorem}

\subsection{Cohomology and local systems}
Assume that $G$ is semisimple, let $\pi:\Gsc\to G$ be its simply connected cover, and put $\Gamma=\ker(\pi)$. For $V\in\RepC(\Gamma)$, let $M_V\in\Invdesc(G)$ be the corresponding connection and let $\mathcal L_V$ be its horizontal local system. 
\begin{theorem}[\cref{thm:derham}]
There is a canonical isomorphism in the derived category of complex vector spaces
\[
  \RGammadR(G,M_V)
  \cong
  \bigl(
    \RGammadR(\Gsc,\cO_{\Gsc})\otimes_{\C}V
  \bigr)^\Gamma.
\]
Consequently, for every $i\geq0$,
\[
  \HdR^i(G,M_V)
  \cong
  \bigl(\HdR^i(\Gsc,\C)\otimes_{\C}V\bigr)^\Gamma
  \cong
  \HdR^i(\Gsc,\C)\otimes_{\C}V^\Gamma.
\]
In particular,
$
  \HdR^0(G,M_V)\cong V^\Gamma.
$
\end{theorem}
\begin{corollary}[\cref{cor:borel}]
Let $d_1,\dots,d_r$ be the fundamental degrees of the Weyl group of $G$, where $r=\operatorname{rank}(G)$. There is a noncanonical isomorphism of graded algebras
\[
  \HdR^*(G,\C)
  \cong
  \Lambda_{\C}(x_1,\dots,x_r),
  \qquad
  \deg(x_j)=2d_j-1.
\]
For $V\in\RepC(\Gamma)$,
\[
  \HdR^*(G,M_V)
 \cong
  \Lambda_{\C}(x_1,\dots,x_r)\otimes_{\C}V^\Gamma
\]
as graded vector spaces.
\end{corollary}
\begin{theorem}[\cref{thm:local-systems}]
The assignment $V\longmapsto\mathcal L_V$ defines an exact tensor equivalence
\[
  \RepC(\Gamma)\isoto\LocSysC(G^{\mathrm{an}}).
\]
Thus every local system has finite monodromy, the category is semisimple, and its simple objects are rank-one local systems indexed by the characters of $\Gamma$.
\end{theorem}
\begin{theorem}[\cref{thm:betti}]
For every $i\geq0$, there is a canonical isomorphism
\[
  H^i(G^{\mathrm{an}},\mathcal L_V)
 \cong
  \bigl(H^i((\Gsc)^{\mathrm{an}},\C)\otimes_{\C}V\bigr)^\Gamma.
\]
\end{theorem}

\subsection*{Relation to previous work}
Classical invariant-connection theory describes invariant principal connections on homogeneous bundles in terms of linear data at a base point~\cite{Wan58}. Here we impose flatness in the algebraic category. Moreover, morphisms are arbitrary algebraic horizontal morphisms, not necessarily equivariant. Thus ``left-invariant'' refers to the existence of an invariant connection form in some algebraic trivialization.

Wei~\cite{Wei25} introduced precisely the free-module formulation corresponding to $\Invtriv(G)$ and obtained classifications for tori, unipotent groups, Borel subgroups of general linear groups, and simply connected semisimple groups. Wei also proposed that pullback from the abelianization should classify all such objects. The present paper shows that this is true for connected reductive groups exactly when the derived subgroup is simply connected, and gives the $\PGL_2$ counterexample otherwise.

On the other hand, Deligne's regular Riemann--Hilbert correspondence identifies regular-singular algebraic connections with complex local systems~\cite[II.5.9]{Del70}, while the fundamental group of a connected complex reductive group is described by the root datum as $X_*(T)/Q^\vee$; see, for example,~\cite{BGA14}. These standard results classify the ambient regular-singular category, but they do not determine which objects admit a left-invariant presentation on a trivial vector bundle.

The contribution of the paper is the bridge between these two levels. The concrete category $\Invtriv(G)$ remembers which monodromy representations are realized on trivial vector bundles. Finite central descent enlarges it to
$
  \Invdesc(G)=\Connrs(G),
$
and passage to direct summands gives
$
  \Kar\bigl(\Invtriv(G)\bigr)\simeq\Connrs(G).
$
Thus left-invariant trivial-bundle algebraic connections form a system of algebraic tensor generators for the entire regular-singular de Rham category of $G$. The derived-monodromy functor measures the remaining obstruction to coming from the abelianization.
\section{Preliminaries and known results}\label{sec:prelim}

\subsection{Algebraic connections, horizontal morphisms, and gauges}
Let $X$ be a smooth affine complex variety, let $\C[X]$ be its coordinate ring and $\Omega_{\C[X]}$ be its module of K\"ahler differentials, and let $D_X$ denote the algebra of algebraic differential operators. For a flat matrix-valued one-form
\[
  \alpha\in M_n(\Omega_{\C[X]}),
  \qquad d\alpha+\alpha\wedge\alpha=0,
\]
write
\[
  M_\alpha=(\cO_X^{\oplus n},d+\alpha)
\]
for the associated vector bundle with integrable connection, viewed equivalently as a locally free left $D_X$-module.

We have the following proposition from
\cite[Proposition~1.1]{Wei25}.

\begin{proposition}[Horizontal morphisms, Wei]
\label{prop:horizontal-morphisms}
Let $\alpha\in M_n(\Omega_{\C[X]})$ and
$\beta\in M_m(\Omega_{\C[X]})$ be flat. Then
$
  \Hom_{D_X}(M_\alpha,M_\beta)
  =
  \left\{
    F\in M_{m,n}(\C[X])
    \ \middle|\
    dF=F\alpha-\beta F
  \right\}.
$
Equivalently, it is enough to impose
$
  v(F)=F\alpha(v)-\beta(v)F
$
for every algebraic vector field $v$ on $X$.
\end{proposition}

\begin{proof}
A matrix $F$ defines an $\cO_X$-linear map between the two trivial vector bundles. Horizontality is the identity
$
  (d+\beta)\circ F=F\circ(d+\alpha),
$
which is equivalent to
$
  dF=F\alpha-\beta F.
$
Evaluating on algebraic vector fields gives the second formulation.
\end{proof}

\begin{remark}[Gauge convention]\label{rem:gauge-convention}
A matrix $F$ represents a morphism $M_\alpha\to M_\beta$ precisely
when
$
  dF=F\alpha-\beta F.
$
If $F$ is invertible, then
$
  \beta=F\alpha F^{-1}-dF\,F^{-1}.
$
All gauge calculations below use this convention.
\end{remark}

Let a linear algebraic group $H$ act on $X$. We write $\Omega_{\C[X]}^H$ for the invariant algebraic one-forms. When $X=H$ and the action is left multiplication, these are the left-invariant forms, and evaluation at the identity gives a canonical isomorphism
$
  \Omega_{\C[H]}^H\isoto \Lie(H)^*, 
$ see \cref{prop:invariant-forms-lie-reps} for more details.

\begin{definition}[Trivial-bundle invariant category]
\label{def:invtriv}
Let $H$ be a complex affine algebraic group. A finite-rank algebraic
connection $M=(E,\nabla)$ on $H$ belongs to $\Invtriv(H)$ if there
exist a finite-dimensional complex vector space $V$ and an algebraic
vector-bundle isomorphism
\[
  \tau:E\xrightarrow{\sim}\cO_H\otimes_{\C}V
\]
such that
$
(\tau\otimes\operatorname{id}_{\Omega_H^1})
  \circ\nabla\circ\tau^{-1}
  =d+\alpha,
$
where
$
  \alpha\in
  \Gamma(H,\Omega_H^1)\otimes_{\C}\End_{\C}(V)
$
is left invariant and satisfies
$
  d\alpha+\alpha\wedge\alpha=0.
$
Necessarily $\dim_{\C}V=\operatorname{rank}(E)$.

The category $\Invtriv(H)$ is the full, replete subcategory of
finite-rank algebraic connections satisfying this condition.
Morphisms are arbitrary algebraic horizontal morphisms.
\end{definition}

\begin{remark}
The vector space $V$ and the trivialization $\tau$ are auxiliary:
neither is part of the data of an object of $\Invtriv(H)$.
\end{remark}

Let $G$ now be a connected reductive group, where $\Gder:=[G,G]$ and $Z(G)^0$ is the connected center of~$G$. Put $Z=Z(G)^0$, and let
\[
  \pi:\widetilde G:=(\Gder)^{\mathrm{sc}}\times Z\longrightarrow G,
  \qquad (s,z)\longmapsto \pi_{\mathrm{sc}}(s)z,
\]
be the canonical finite central isogeny.

\begin{definition}[Canonical-cover descent category]\label{def:invdesc}
A finite-rank algebraic connection $M$ on $G$ is \emph{descent-invariant} if
$
  \pi^*M\in\Invtriv(\widetilde G).
$
The full subcategory of such connections is denoted $\Invdesc(G)$.
\end{definition}

\begin{remark}
In $\Invdesc(G)$, only the pullback of the underlying vector bundle is
required to be algebraically trivial. The vector bundle on $G$ itself
may be nontrivial.
\end{remark}

\begin{remark}\label{rem:invdesc-flat}
Flatness need not be imposed separately in
\cref{def:invdesc}. Indeed, if $M=(E,\nabla)$ and
$
  F_\nabla=\nabla^2:E\longrightarrow E\otimes\Omega_G^2
$
is its curvature, then
$
  F_{\pi^*\nabla}=\pi^*F_\nabla.
$
 Since $\pi^*M\in\Invtriv(\widetilde G)$, the connection
$\pi^*\nabla$ is flat, and hence $\pi^*F_\nabla=0$.
The morphism $\pi$ is finite \'{e}tale and surjective, hence
faithfully flat. Therefore pullback is faithful on morphisms of
vector bundles, and consequently $F_\nabla=0$.
\end{remark}

Throughout, when a finite central kernel such as $\Gamma$ or $K$
appears in a representation category or in a fundamental-group
sequence, we identify it with its finite group of complex points.

\subsection{Relation with Wei's notation}\label{subsec:wei-notation}
Wei writes $D(X)$ for the algebra denoted here by $D_X$. Given a flat matrix-valued one-form
$
  \alpha\in M_n(\Omega_{\C[X]}),
  \; d\alpha+\alpha\wedge\alpha=0,
$
the module $M_\alpha$ is precisely the connection
$
  (\cO_X^{\oplus n},d+\alpha).
$

When $X=G$ and $\alpha$ is left invariant, Wei's invariant algebraic $D$-modules of rank $n$ are therefore exactly the objects of $\Invtriv_n(G)$. In particular,
\[
  \frac{\text{{Invariant $D$-modules of rank $n$ on $G$} in the sense of~\cite{Wei25}}}
       {\text{$D_G$-module isomorphism}}
 =
  [\Invtriv_n(G)].
\]
Passing to the full replete category $\Invtriv(G)$ does not change the set of isomorphism classes. The category $\Invdesc(G)$ has no counterpart in \cite{Wei25}: its objects may have nontrivial underlying vector bundles and arise by finite central descent.

For reference, Wei's conjecture may be written in the present notation as follows.
\begin{conjecture}[Wei]\label{conj:wei}
Let $G$ be a connected linear algebraic group, and let $q:G\twoheadrightarrow G/[G,G]$ be the quotient map. For every $n\geq1$, pullback induces a bijection
\[
  q^*: [\Invtriv_n(G/[G,G])]
  \longrightarrow
  [\Invtriv_n(G)].
\]
\end{conjecture}
The conjecture concerns only the trivial-bundle category. The $\PGL_2$ example in \cref{sec:semisimple} gives a counterexample, while \cref{cor:wei-sc-derived} identifies an important class of reductive groups for which the conjecture holds. \cref{cor:wei-iff} gives a sharper result for this conjecture.

When the canonical cover is an isomorphism, for example when $G$ is a torus or a simply connected semisimple group, the two categories agree.

We recall the following proposition; see \cite[Proposition 1.3]{Wei25}.
\begin{proposition}[Wei]\label{prop:invariant-forms-lie-reps}
For every $n \geq 1$, evaluation at the identity gives a canonical bijection between
\[
  \left\{
  \alpha\in M_n\bigl(\Omega_{\C[H]}^H\bigr)
  \ \middle|\ 
  d\alpha+\alpha\wedge\alpha=0
  \right\}
\]
and the set of Lie algebra homomorphisms
$
  \rho:\Lie(H)\longrightarrow\mathfrak{gl}_n.
$
\end{proposition}

\begin{proof}
A left-invariant form is determined by its value $\rho$ at the identity. For left-invariant vector fields $u,v$,
$
  (d\alpha+\alpha\wedge\alpha)(u,v)
  =-\rho([u,v])+[\rho(u),\rho(v)].
$
Thus flatness is equivalent to the homomorphism condition.
\end{proof}

For a Lie algebra representation $\rho:\Lie(H)\to\mathfrak{gl}(V)$, write $M_\rho$ for the corresponding object of $\Invtriv(H)$.

We have the following theorem from \cite[Theorem~6.1]{BF18}.
\begin{theorem}[Billig-Futorny]\label{thm:BF}
For a complex affine algebraic group $H$ with Lie algebra $\mathfrak h$,
\[
  \Der_{\C}(\C[H])\cong\C[H]\otimes_{\C}\mathfrak h,
\]
where $\mathfrak h$ is identified with the left-invariant algebraic vector fields.
\end{theorem}
Consequently, in \cref{prop:horizontal-morphisms} it is enough to test the horizontal equation on $v\in\mathfrak h$. Wei proved the following proposition in \cite{Wei25} for the special case of weakly exponential algebraic groups. Here we state and prove a version for arbitrary complex affine algebraic groups.
\begin{proposition}\label{prop:integrable-gauge}
Let $H$ be a complex affine algebraic group and let $\rho:\Lie(H)\to\mathfrak{gl}(V)$ be a Lie algebra homomorphism. Suppose that $\rho$ is the differential of an algebraic group homomorphism
$
  \Phi:H\longrightarrow\GL(V).
$
Then $M_\rho$ is algebraically isomorphic to the trivial connection.
\end{proposition}

\begin{proof}
Define $X(h)=\Phi(h)$. If $v\in\Lie(H)$ is viewed as a left-invariant vector field, then
\[
  v(X)(h)
  =\left.\frac{d}{dt}\right|_{t=0}\Phi(h\exp(tv))
  =X(h)\rho(v).
\]
Thus $dX=X\alpha_\rho$, and \cref{rem:gauge-convention} shows that $X$ is an isomorphism $M_\rho\to M_0$.
\end{proof}

\subsection{The Maurer--Cartan form}

\begin{definition}
Let $H$ be a smooth algebraic group with Lie algebra $\mathfrak h$. The left Maurer--Cartan form is the unique left-invariant $\mathfrak h$-valued one-form $\theta_H$ satisfying
$
  (\theta_H)_e=\operatorname{id}_{\mathfrak h}.
$
Equivalently,
\[
  (\theta_H)_h=(dL_{h^{-1}})_h:T_hH\longrightarrow T_eH=\mathfrak h.
\]
\end{definition}

\begin{proposition}\label{prop:MC}
The Maurer--Cartan form satisfies
$
  d\theta_H+\theta_H\wedge\theta_H=0.
$
Every left-invariant flat algebraic connection on a trivial vector bundle has the form
\[
  d+(\rho\otimes{\operatorname{id}_{\Omega_H^1}})(\theta_H)\]
for a unique Lie algebra homomorphism $\rho:\mathfrak h\to\mathfrak{gl}(V)$.
\end{proposition}

\begin{proof}
The first part is well known; see, for example, \cite[Theorem 15.3]{Tu}. We include a proof for completeness. For left-invariant vector fields $u,v$, the functions
$\theta_H(u)$ and $\theta_H(v)$ are constant. Hence
$
  d\theta_H(u,v)
  =-\theta_H([u,v])
  =-[u,v].
$
With the convention
$
  (\theta_H\wedge\theta_H)(u,v)
  =
  [\theta_H(u),\theta_H(v)],
$
the second term equals $[u,v]$, proving the Maurer--Cartan equation.

The second assertion follows because every left-invariant
$\End(V)$-valued one-form is uniquely of the form
$(\rho\otimes{\operatorname{id}_{\Omega_H^1}})(\theta_H)$, and by
\cref{prop:invariant-forms-lie-reps} its flatness is equivalent to
$\rho$ being a Lie algebra homomorphism.
\end{proof}
\begin{remark}
For $G=\GL_n$, where $g$ denotes the tautological matrix-valued
coordinate, the form $g^{-1}dg$ is the left Maurer--Cartan form.
\end{remark}
\subsection{Tori and comparison with previous work}
We first recall the notion of a regular-singular algebraic connection.
\begin{definition}[Regular-singular connection]\label{def:regular-singular}
Let $X$ be a smooth complex algebraic variety. A finite-rank flat algebraic connection $M$ on $X$ is \emph{regular-singular} if, for every morphism $f:C\to X$ from a smooth connected algebraic curve, the pullback $f^*M$ has regular singularities along the boundary of a smooth projective compactification of $C$.

We denote by $\Connrs(X)$ the category of finite-rank regular-singular algebraic connections on~$X$.
\end{definition}
\begin{remark}
By Deligne's curve and compactification criteria, the above definition of regular-singular is equivalent to his intrinsic notion of regularity and is independent of the chosen compactification~\cite[II.4.4--II.4.5]{Del70}. 
\end{remark}

The proofs below use only two elementary inputs: the classification on a torus and the gauge triviality on a simply connected semisimple group. We prove both directly. Wei also treats unipotent groups and Borel subgroups in \cite{Wei25}.

The following theorem recovers the torus classification of \cite[Theorem~0.3]{Wei25}.
\begin{theorem}[Tori]\label{thm:wei-torus}
Let $T\cong(\C^*)^r$ be an algebraic torus. For every $n \geq 1$, monodromy induces a natural bijection
\[
  [\Invtriv_n(T)]=
  \Hom\bigl(\pi_1(T^{\mathrm{an}},e),\GL_n\bigr)/\GL_n.
\]
Equivalently, every finite-rank complex local system on $T^{\mathrm{an}}$ arises as the horizontal local system of an invariant trivial-bundle algebraic connection.
\end{theorem}

\begin{proof}
\emph{1. From connections to monodromy.} Choose coordinates $z_1,\dots,z_r$ on $T$. A left-invariant flat algebraic connection on a trivial vector bundle has the form
\[ d+\sum_{i=1}^r B_i\frac{dz_i}{z_i},
  \qquad [B_i,B_j]=0,
\]
and extends to a logarithmic connection on $(\mathbb P^1)^r$ along its boundary. In particular, it is regular-singular.

\emph{2. Surjectivity: from representations to connections.} Conversely, let a representation of $\pi_1(T^{\mathrm{an}},e)\cong~\mathbb Z^r$ be given by commuting matrices 
\[A_1,\dots,A_r\in\GL(V).\] 
The commutative algebra 
\[A=\C[A_1,\dots,A_r]\subset\End(V)\]
is finite-dimensional. Each $A_i$ is a unit of $A$: by the Cayley--Hamilton theorem, $A_i^{-1}$ is a polynomial in $A_i$. We write the Artinian decomposition as follows: 
\[
A\cong\prod_{\nu=1}^m A_\nu\text{, where each $A_\nu$ is a local Artinian $\C$-algebra with maximal ideal $\mathfrak m_\nu$}.\]

If $u=(u_\nu)_\nu\in A^*$, then for each $\nu$ there are $\lambda_\nu\in\C^*$ and a nilpotent element $N_\nu\in\mathfrak m_\nu$ such that $u_\nu=\lambda_\nu(1+N_\nu)$.
Choose $c_\nu\in\C$ with $\exp(c_\nu)=\lambda_\nu$ and set
\[
L_\nu:=c_\nu+\sum_{k\geq1}\frac{(-1)^{k+1}}{k}N_\nu^k.
\]
The sum is finite and $\exp(L_\nu)=u_\nu$. Reassembling the $L_\nu$ gives $L\in A$ with $\exp(L)=u$. Applying this to the matrices $A_i$ gives commuting elements $L_i\in A$ with $\exp(L_i)=A_i$.

The invariant connection
\[
  d-\frac{1}{2\pi i}\sum_{i=1}^r L_i\frac{dz_i}{z_i}\]
has monodromy $A_i$ around the $i$-th standard loop. This proves surjectivity.

\emph{3. Injectivity: from monodromy to algebraic isomorphisms.} If two such invariant connections have conjugate monodromy representations, then their associated local systems are isomorphic. Since both connections are regular-singular, the full faithfulness of Deligne's regular-singular Riemann--Hilbert correspondence implies that this isomorphism is induced by an algebraic isomorphism of connections. This proves injectivity.
\end{proof}

\begin{lemma}[Finite-isogeny descent on tori]\label{lem:torus-descent}
Let $q:T'\to T$ be a finite isogeny of complex algebraic tori. If an object of $\Invtriv(T')$ is equipped with descent data for $q$, then its descended flat bundle on $T$ is isomorphic to an object of $\Invtriv(T)$.
\end{lemma}

\begin{proof}
Let $M'\in\Invtriv(T')$, and let $M$ be its descent to $T$, so that
$q^*M\cong M'$. Since an invariant algebraic connection on a torus
extends logarithmically to a smooth toric compactification, $M'$ is
regular-singular. Regular singularity is reflected by finite étale
surjective pullback~\cite[II.4.6(iii)]{Del70}; hence $M$ is regular-singular.

By \cref{thm:wei-torus}, there exists $N\in\Invtriv(T)$ whose
monodromy is isomorphic to that of~$M$. Both $M$ and $N$ are regular-singular, so the full faithfulness of Deligne's regular-singular
Riemann--Hilbert correspondence yields an algebraic isomorphism
$M\cong N$~\cite[II.5.9]{Del70}.
\end{proof}

The following theorem gives a categorical formulation of Wei's result \cite[Theorem~0.5]{Wei25}.
\begin{theorem}[Simply connected semisimple groups]\label{thm:wei-sc}
If $S$ is a connected, simply connected, semisimple group, then for every
$n\geq1$, every rank-$n$ object of $\Invtriv(S)$ is isomorphic to the
rank-$n$ trivial connection.
\end{theorem}

\begin{proof}
An object of $\Invtriv(S)$ is represented by a Lie algebra homomorphism $\rho:\Lie(S)\to\mathfrak{gl}(V)$. Since $S$ is simply connected and semisimple, $\rho$ integrates uniquely to an algebraic representation $\Phi:S\to\GL(V)$. The gauge $X(s)=\Phi(s)$ trivializes the connection by \cref{prop:integrable-gauge}.
\end{proof}

\subsection{The descent category and regular-singular connections}
We denote by $\LocSysC(G^{\mathrm{an}})$ the category of finite-rank complex local systems.

\begin{theorem}\label{thm:rs}
Let $G$ be a connected complex reductive algebraic group. Then
\[
  \Invdesc(G)=\Connrs(G)
\]
as full subcategories of the category of finite-rank algebraic connections. Consequently, horizontal sections and monodromy induce exact tensor equivalences
\[
  \Invdesc(G)\isoto \LocSysC(G^{\mathrm{an}})
  \isoto \RepC\bigl(\pi_1(G^{\mathrm{an}},e)\bigr).
\]
\end{theorem}

\begin{proof}
Put $S=(\Gder)^{\mathrm{sc}}$, $Z=Z(G)^0$, and let
$
  \pi:S\times Z\longrightarrow G
$
be the canonical finite central isogeny.

Let $M\in\Invdesc(G)$. By definition, $\pi^*M$ is isomorphic to the invariant connection associated with a Lie algebra homomorphism
\[
  \rho:\Lie(S)\oplus\Lie(Z)\longrightarrow\mathfrak{gl}(V).\]
Write $\rho_S$ and $\rho_Z$ for the restrictions. Their images commute. Since $S$ is simply connected and semisimple, $\rho_S$ integrates uniquely to an algebraic representation
$
  \Phi:S\longrightarrow\GL(V).
$
The gauge $X(s,z)=\Phi(s)$ removes the $S$-component. Since $\rho_Z(\Lie(Z))$ commutes with the connected subgroup $\Phi(S)$, the $Z$-component is unchanged. Hence
$
  \pi^*M\cong\pr_Z^*N_Z
$
for some $N_Z\in\Invtriv(Z)$.

The torus connection $N_Z$ is regular-singular, and regular singularity is preserved by pullback, so
$\pi^*M$ is regular-singular. Since $\pi$ is finite étale and
surjective, regular singularity is reflected by pullback along
$\pi$~\cite[II.4.6(iii)]{Del70}. Therefore $M$ is regular-singular.

Conversely, let $M\in\Connrs(G)$. Then $\pi^*M$ is regular-singular. The complex manifold $S^{\mathrm{an}}$ is simply connected~\cite[Appendix D, Proposition D.4.1]{Con14}, so projection
$
  \pr_Z:S^{\mathrm{an}}\times Z^{\mathrm{an}}\longrightarrow Z^{\mathrm{an}}
$
induces an isomorphism on fundamental groups. More explicitly, restriction along
\[
  \iota_Z:Z^{\mathrm{an}}\hookrightarrow
  S^{\mathrm{an}}\times Z^{\mathrm{an}},\qquad z\mapsto(e,z),
\]
is inverse to pullback along $\pr_Z$ at the level of local systems.
Hence the local system underlying $\pi^*M$ is isomorphic to
$\pr_Z^*L_Z$ for some local system $L_Z$ on $Z^{\mathrm{an}}$.

By \cref{thm:wei-torus}, there exists
$N_Z\in\Invtriv(Z)$ whose horizontal local system is $L_Z$.
Consequently, the horizontal local system of $\pr_Z^*N_Z$ is
$\pr_Z^*L_Z$, which is isomorphic to that of $\pi^*M$. Both
$\pi^*M$ and $\pr_Z^*N_Z$ are regular-singular, so the full
faithfulness of Deligne's regular-singular Riemann--Hilbert
correspondence gives an algebraic isomorphism
\[
  \pi^*M\cong\pr_Z^*N_Z.
\]
Since $\pr_Z^*N_Z\in\Invtriv(S\times Z)$, it follows that
$M\in\Invdesc(G)$.

The remaining equivalences follow from Deligne's regular-singular
Riemann--Hilbert correspondence~\cite[II.5.9]{Del70} and the
monodromy equivalence between complex local systems on
$G^{\mathrm{an}}$ and finite-dimensional representations of
$\pi_1(G^{\mathrm{an}},e)$. Both equivalences are exact and
compatible with tensor products.
\end{proof}

\begin{corollary}\label{cor:translation-invariance}
For every $M\in\Invdesc(G)$ and every $g\in G(\C)$, there exists an algebraic isomorphism
$
  L_g^*M\cong M.
$
\end{corollary}

\begin{proof}
Because $G^{\mathrm{an}}$ is connected, $L_g$ is homotopic to the identity. It therefore preserves every local system up to isomorphism, and \cref{thm:rs} gives the algebraic statement.
\end{proof}

\begin{remark}\label{rem:nonequivariant}
The isomorphism in \cref{cor:translation-invariance} is generally noncanonical. The statement does not supply coherent isomorphisms for varying $g$, and it should not be confused with a $G$-equivariant structure in the usual theory of equivariant $D$-modules.
\end{remark}

\section{The semisimple case}\label{sec:semisimple}

Let $G$ be a connected semisimple complex algebraic group, let
$
  \pi:\Gsc\twoheadrightarrow G
$
be its simply connected central cover, and put
$
  \Gamma:=\ker(\pi)\subset Z(\Gsc).
$
The group $\Gamma$ is finite and abelian.

\subsection{The trivial-bundle classification}

\begin{proposition}\label{prop:semisimple-triv}
Fix $n\geq1$. Restriction to $\Gamma$ induces a natural bijection
\[
  [\Invtriv_n(G)]
  =
  \operatorname{im}\!\left(
  \frac{\Hom_{\mathrm{alg.grp}}(\Gsc,\GL_n)}{\GL_n}
  \longrightarrow
  \frac{\Hom(\Gamma,\GL_n)}{\GL_n}
  \right).
\]
\end{proposition}

\begin{proof}
An object of $\Invtriv_n(G)$ is represented by a Lie algebra homomorphism
\[\rho:\g\longrightarrow\mathfrak{gl}_n.\]
Using the isomorphism $d\pi:\Lie(\Gsc)\xrightarrow{\sim}\g$, we
regard $\rho$ as a representation of $\Lie(\Gsc)$. Since $\Gsc$ is
simply connected and semisimple, it integrates uniquely to an
algebraic representation
\[
  \Phi:\Gsc\longrightarrow\GL_n.
\]
We associate to $M_\rho$ the conjugacy class of $\Phi|_\Gamma$.

Let $M_\rho$ and $M_\sigma$ correspond to $\Phi$ and $\Psi$, respectively. After pullback to $\Gsc$, the gauge transformations $\Phi$ and $\Psi$ trivialize the corresponding connections. Hence every horizontal isomorphism between the pullbacks has the form
\[
  \widetilde Y(g)=\Psi(g)^{-1}B\Phi(g),
  \qquad B\in\GL_n(\C).
\]
It descends to $G=\Gsc/\Gamma$ if and only if 
\[\widetilde Y(g\gamma)=\widetilde Y(g)\qquad(\gamma\in\Gamma).\] Since $\Gamma$ is central,
this condition is equivalent to
\[
  B\Phi(\gamma)=\Psi(\gamma)B
  \qquad(\gamma\in\Gamma).\]
Thus $M_\rho\cong M_\sigma$ if and only if the restrictions to $\Gamma$ are conjugate.

Finally, every element of the displayed image is represented by
$\Phi|_\Gamma$ for some algebraic representation
$\Phi:\Gsc\to\GL_n$, and $d\Phi$ defines an object of
$\Invtriv_n(G)$. This proves the claimed bijection.
\end{proof}

\begin{remark}\label{rem:semisimple-distinction}
Only representations of $\Gamma$ that occur as restrictions of algebraic representations of $\Gsc$ appear in $\Invtriv(G)$. The finite-descent category allows every finite-dimensional representation of $\Gamma$.
\end{remark}

\subsection{Finite descent and the category \texorpdfstring{$\RepC(\Gamma)$}{Rep(Gamma)}}

\begin{lemma}[Integration]\label{lem:integration}
Let $\Gsc$ be a simply connected semisimple complex algebraic group with Lie algebra $\mathfrak{g}$. Every Lie algebra homomorphism
$
  \Lie(\Gsc)\longrightarrow\mathfrak{gl}(V)
$
is the differential of a unique algebraic group homomorphism $\Gsc\to\GL(V)$.
\end{lemma}

\begin{proof}
This is the standard integration theorem for simply connected semisimple algebraic groups; see Serre~\cite[p.~152, Theorem~1]{Ser92}.
\end{proof}

\begin{lemma}[Automorphisms of the trivial connection]\label{lem:aut-trivial}
Let $X$ be a connected smooth complex algebraic variety. Then
\[ \Aut_{D_X}(\cO_X\otimes_{\C}V,d)=\GL(V).\]
More generally, every horizontal morphism between two trivial connections is a constant linear map.
\end{lemma}

\begin{proof}
A horizontal morphism is represented by a matrix $F$ of regular functions satisfying $dF=0$. On a connected smooth complex variety, every regular function with vanishing differential is constant.
\end{proof}

\begin{lemma}[Finite \'{e}tale Galois descent]\label{lem:etale-descent}
Let $p:Y\to X$ be a finite \'{e}tale Galois cover of smooth complex varieties with group $\Delta$. Pullback induces an equivalence between algebraic connections on $X$ and $\Delta$-equivariant algebraic connections on $Y$.
\end{lemma}

\begin{proof}
Since $p$ is finite \'etale and surjective, it is faithfully flat.
Hence vector bundles, their morphisms, and descent data satisfy
effective faithfully flat descent. For a finite Galois cover,
descent data are equivalent to compatible $\Delta$-equivariant
structures.

Let $(E_Y,\nabla_Y)$ be a $\Delta$-equivariant algebraic connection
on $Y$. The underlying $\Delta$-equivariant vector bundle descends
to a vector bundle $E$ on $X$, with $p^*E\cong E_Y$. Since $p$ is \'etale, there is a canonical identification $\Omega_Y^1\cong p^*\Omega_X^1$.
Thus the $\Delta$-equivariant morphism
\[
  \nabla_Y:E_Y\longrightarrow E_Y\otimes\Omega_Y^1
\]
descends uniquely to a morphism
\[
  \nabla:E\longrightarrow E\otimes\Omega_X^1.
\]
The Leibniz identity for $\nabla$ can be checked after the faithfully
flat pullback $p^*$, where it becomes the Leibniz identity for
$\nabla_Y$. Hence $\nabla$ is an algebraic connection.

Moreover, under the identification
$\Omega_Y^2\cong p^*\Omega_X^2$, one has
\[
  F_{\nabla_Y}=p^*F_\nabla.
\]
Since faithfully flat pullback detects the vanishing of morphisms,
$\nabla_Y$ is flat if and only if $\nabla$ is flat. Finally,
morphisms of vector bundles descend uniquely, and horizontality can
also be checked after pullback. This proves full faithfulness and
essential surjectivity.
\end{proof}

\begin{theorem}\label{thm:semisimple-desc}
For every $n\geq1$, there is a natural bijection
\[
  [\Invdesc_n(G)]
  =
  \Hom(\Gamma,\GL_n)/\GL_n.
\]
Moreover, there is an exact tensor equivalence
$
  \Invdesc(G)\isoto\RepC(\Gamma).
$
\end{theorem}

\begin{proof}
Let $M\in\Invdesc(G)$. Its pullback to $\Gsc$ lies in $\Invtriv(\Gsc)$ and is therefore trivial by \cref{thm:wei-sc}.
By \cref{lem:etale-descent}, the object $M$ is equivalent to a $\Gamma$-linearization of
$
  (\cO_{\Gsc}\otimes V,d).
$ We use the right deck transformations $t_\gamma(g)=g\gamma^{-1}$.
For each $\gamma\in\Gamma$, the linearization map 
\[
A_\gamma:t_\gamma^*(\cO_{\Gsc}\otimes V,d)
  \longrightarrow(\cO_{\Gsc}\otimes V,d)
\]
is a horizontal automorphism of the trivial connection and hence, by \cref{lem:aut-trivial}, a constant matrix $A_\gamma\in\GL(V)$. 
The group $\Gamma$ is abelian. Interchanging $\gamma_1$ and $\gamma_2$ therefore gives $A_{\gamma_1\gamma_2}=A_{\gamma_1}A_{\gamma_2}$, so $\gamma\mapsto A_\gamma$ is an ordinary representation. 
Changing the trivialization conjugates this representation. This proves the bijection on isomorphism classes.

Evaluation at the identity identifies trivial connections on $\Gsc$ with finite-dimensional vector spaces and is compatible with tensor products.
Passing to $\Gamma$-equivariant objects gives the exact tensor equivalence with $\RepC(\Gamma)$.
\end{proof}

\subsection{The example of \texorpdfstring{$\PGL_2$}{PGL2} and Wei's conjecture}

Let
\[
  G=\PGL_2,
  \qquad
  \Gsc=\SL_2,
  \qquad
  \Gamma=\mu_2=\{\pm I_2\}.
\]
\subsubsection{The rank-one case}

Every algebraic character $\SL_2\to\Gm$ is trivial. Hence \cref{prop:semisimple-triv} shows that $\Invtriv_1(\PGL_2)$ has only the trivial isomorphism class. By \cref{thm:semisimple-desc}, however, $\Invdesc_1(\PGL_2)$ has two classes, corresponding to the trivial and sign characters of $\mu_2$.

\begin{example}[The descended sign line]\label{ex:sign-line}
Let $\C_{\mathrm{sgn}}$ be the one-dimensional sign representation of $\mu_2$. The associated bundle
$
  \mathcal L_{\mathrm{sgn}}
  :=\SL_2\times^{\mu_2}\C_{\mathrm{sgn}}
  \longrightarrow\PGL_2
$
is equipped with the connection descended from the trivial connection on $\SL_2\times\C$. Its monodromy is the nontrivial character of $\pi_1(\PGL_2,e)\cong\mu_2$. This object belongs to $\Invdesc(\PGL_2)$ but not to $\Invtriv(\PGL_2)$. Its underlying line bundle is nontrivial: a trivialization would pull back to a unit $f\in\cO(\SL_2)^*$ satisfying $f(-g)=-f(g)$, whereas every unit on the connected semisimple group $\SL_2$ is constant by the following lemma.
\end{example}
\begin{lemma}\label{lem:semisimple-units}
If $S$ is a connected semisimple complex algebraic group, then
\[
  \cO(S)^*\cong\C^*.
\]
\end{lemma}

\begin{proof}
Rosenlicht's unit theorem (\cite{Ros56}) identifies $\cO(S)^*/\C^*$ with the algebraic character group $X^*(S)$. A connected semisimple group has no nontrivial algebraic characters.
\end{proof}

\subsubsection{The rank-two case}
By complete reducibility and the classification of irreducible representations of $\SL_2$, every two-dimensional algebraic representation is either the direct sum of two trivial representations or the standard representation. Thus the central element $-I_2$ acts either as $I_2$ or as $-I_2$, and the mixed action
$
  -I_2\longmapsto \diag(1,-1)
$
appears only in $\Invdesc(\PGL_2)$.

\begin{example}[Descended from the standard representation]\label{ex:std}
Let $\Phi:\SL_2\hookrightarrow\GL_2$ be the standard representation and let
\[\rho=d\Phi:\mathfrak{pgl}_2\cong\mathfrak{sl}_2\longrightarrow\mathfrak{gl}_2.\]
The differential defines a left-invariant algebraic connection $M_\rho$ on the trivial rank-two vector bundle over $\PGL_2$. Since $\Phi(-I_2)=-I_2$, \cref{prop:semisimple-triv} gives $M_\rho\not\cong M_0$. Because $\PGL_2^{\mathrm{ab}}=\{1\}$, this is a counterexample to \cref{conj:wei} in Wei's original trivial-bundle category. Thus the source of the pullback map in \cref{conj:wei} has only the trivial rank-two class, while its target also contains the nontrivial class of $M_\rho$.
\end{example}
\subsubsection{Splitting in the descent category}
Under the equivalence $\Invdesc(\PGL_2)\simeq\RepC(\mu_2)$, the object $M_\rho$ corresponds to $\mathrm{sgn}^{\oplus2}$. Consequently,
$
  M_\rho\cong\mathcal L_{\mathrm{sgn}}^{\oplus2}
  \;\text{in }\Invdesc(\PGL_2).
$
The projector onto either sign summand is therefore an idempotent endomorphism of $M_\rho$ that does not split inside $\Invtriv(\PGL_2)$. This gives the basic model for the Karoubi-completion theorem below.

\subsection{Karoubi completion and internal Ext groups}
\begin{definition}
A category is \emph{idempotent complete} if every idempotent endomorphism splits. Its \emph{idempotent completion}, or \emph{Karoubi envelope}, is obtained by formally adjoining the images of all idempotent endomorphisms.
\end{definition}
\begin{corollary}\label{cor:semisimple-kar}
The inclusion $\Invtriv(G)\hookrightarrow\Invdesc(G)$ induces a tensor equivalence
\[
  \Kar\bigl(\Invtriv(G)\bigr)\isoto\Invdesc(G)\simeq\RepC(\Gamma).
\]
\end{corollary}

\begin{proof}
Under the equivalence of \cref{thm:semisimple-desc}, $\Invtriv(G)$ identifies with the full subcategory of $\RepC(\Gamma)$ consisting of restrictions of algebraic representations of $\Gsc$.

Because $\Gamma$ is finite abelian and $\C$ has characteristic zero, every irreducible representation of $\Gamma$ is a character $\chi:\Gamma\to\C^*$. The inclusion $\Gamma\subset Z(\Gsc)$ induces a surjection
\[ X^*\bigl(Z(\Gsc)\bigr)\twoheadrightarrow X^*(\Gamma).
\]
Choose a character of $Z(\Gsc)$ whose restriction to $\Gamma$ is $\chi$. Fix a maximal torus of $\Gsc$, and let $P$ and $Q$ denote its weight
and root lattices. Under the standard identification
\[
  X^*\bigl(Z(\Gsc)\bigr)\cong P/Q,\]
choose a weight $\lambda\in P$ representing this central character. Put
\[
2\rho:=\sum_{\alpha\in\Phi^+}\alpha.
\]
Then $2\rho\in Q$ and $\langle2\rho,\alpha_i^\vee\rangle=2$
for every simple coroot $\alpha_i^\vee$. Consequently, $\lambda+N(2\rho)$ is dominant for $N\gg0$ and has
the same class as $\lambda$ in $P/Q$.

Let $W_\chi$ be the irreducible highest-weight representation with this dominant highest weight. 
Its restriction to $\Gamma$ is
\[
  W_\chi|_\Gamma\cong\chi^{\oplus\dim W_\chi}.
\]
In particular, $\chi$ is a direct summand of the restriction of an
algebraic representation of $\Gsc$. Every finite-dimensional representation of $\Gamma$ is a finite
direct sum of characters. Taking the corresponding direct sums of
the representations $W_\chi$, we conclude that every object of
$\RepC(\Gamma)$ is a direct summand of an object belonging to the image of
$\Invtriv(G)$.

Since $\RepC(\Gamma)$ is idempotent complete, the
universal property of the Karoubi completion gives a tensor
equivalence
\[
  \Kar\bigl(\Invtriv(G)\bigr)\simeq\RepC(\Gamma).
\]
Combining this with \cref{thm:semisimple-desc} proves the result.
\end{proof}

\begin{corollary}\label{cor:semisimple-ext}
The abelian category $\Invdesc(G)$ is semisimple, and for all $M,N\in\Invdesc(G)$,
\[
  \Ext^i_{\Invdesc(G)}(M,N)=0
  \qquad(i>0).
\]
\end{corollary}

\begin{proof}
By \cref{thm:semisimple-desc}, there is an exact equivalence
\[
  \Invdesc(G)\simeq\RepC(\Gamma).
\]
Since $\Gamma$ is finite and $\C$ has characteristic zero,
Maschke's theorem implies that $\RepC(\Gamma)$ is semisimple.
Therefore $\Invdesc(G)$ is a semisimple abelian category. In
particular, every object is projective and injective, and hence
\[
  \Ext^i_{\Invdesc(G)}(M,N)=0
  \qquad(i>0).\qedhere
\]
\end{proof}

\begin{remark}\label{rem:internal-ext}
The groups in \cref{cor:semisimple-ext} are internal Ext groups in
the abelian category $\Invdesc(G)$. They need not agree with Ext
groups in the ambient category of algebraic $D_G$-modules.
\end{remark}

\section{An illustrative example: the general linear group case}\label{sec:glr}

Let $r\geq1$ and let
$
  \Det:\GL_r\longrightarrow\Gm
$
be the determinant. Since $\pi_1(\GL_r(\C),e)\cong\mathbb Z\cong\pi_1(\C^*,e)$ and the determinant induces an isomorphism on fundamental groups, \cref{thm:rs} immediately gives the descent-category statement below. We also give an explicit algebraic proof as it explains why the equivalence restricts to the trivial-bundle category.

\begin{theorem}\label{thm:det-equivalence}
Determinant pullback induces exact tensor equivalences
\[
  \Det^*: \Invdesc(\Gm)\isoto\Invdesc(\GL_r),
  \qquad
  \Det^*: \Invtriv(\Gm)\isoto\Invtriv(\GL_r).
\]
Consequently,
$
  [\Invdesc_n(\GL_r)]
  =
  \Hom(\mathbb Z,\GL_n)/\GL_n.
$
\end{theorem}
The proof occupies the next three subsections.
\subsection{The central cover and the semisimple gauge}

Consider the finite central isogeny
\[
  p:\Gm\times\SL_r\longrightarrow\GL_r,
  \qquad
  (t,s)\longmapsto ts.
\]
Its kernel is
$
  \mu_r=
  \{(\zeta,\zeta^{-1}I_r)\mid \zeta^r=1\}.
$

Let
\[
  \theta_{\GL_r}=g^{-1}dg,
  \qquad
  \omega=t^{-1}dt,
  \qquad
  \theta_{\SL_r}=s^{-1}ds
\]
denote the left Maurer--Cartan forms on $\GL_r$, $\Gm$, and
$\SL_r$, respectively.  We suppress the pullback notation for
$\omega$ and $\theta_{\SL_r}$ on $\Gm\times\SL_r$.

Since $p(t,s)=tI_r\,s$, we have
\begin{align*}
  p^*\theta_{\GL_r}
  &=(ts)^{-1}d(ts)\\
  &=s^{-1}t^{-1}\bigl((dt)s+t\,ds\bigr)\\
  &=I_r\,\frac{dt}{t}+s^{-1}ds\\
  &=I_r\,\omega+\theta_{\SL_r}.
\end{align*}

Let
$
  \rho:\mathfrak{gl}_r\longrightarrow\mathfrak{gl}(V)
$
be a Lie algebra representation, and let
$
  A:=\rho(I_r),
  \;
  \rho_{\mathrm{ss}}:=\rho|_{\mathfrak{sl}_r}.
$
The corresponding invariant connection form on $\GL_r$ is
$
  \alpha=\rho(\theta_{\GL_r}).
$
Therefore
\begin{equation}\label{eq:pullback-glr-form}
  p^*\alpha
  =
  A\omega+\rho_{\mathrm{ss}}(\theta_{\SL_r}).
\end{equation}

Because $I_r$ is central in $\mathfrak{gl}_r$, for every
$x\in\mathfrak{sl}_r$ we have
\[
  [A,\rho_{\mathrm{ss}}(x)]
  =
  [\rho(I_r),\rho(x)]
  =
  \rho([I_r,x])
  =0.
\]
Thus $A$ commutes with the image of $\rho_{\mathrm{ss}}$.

Since $\SL_r$ is simply connected and semisimple,
$\rho_{\mathrm{ss}}$ integrates uniquely to an algebraic
representation
$
  \Phi:\SL_r\longrightarrow\GL(V).
$
As $\Phi(\SL_r)$ is connected, the preceding infinitesimal
commutation relation implies that $A$ commutes with
$\Phi(\SL_r)$.

Define
$
  X_0:\Gm\times\SL_r\longrightarrow\GL(V),
  \;
  X_0(t,s)=\Phi(s).
$
Then, we have
$
  X_0^{-1}dX_0
  =
  \rho_{\mathrm{ss}}(\theta_{\SL_r}).
$ Equivalently, this implies 
$
  dX_0
  =
  X_0\,\rho_{\mathrm{ss}}(\theta_{\SL_r}).
$
Using the gauge convention
\[
  \beta=X_0(p^*\alpha)X_0^{-1}-dX_0\,X_0^{-1},
\]
together with \eqref{eq:pullback-glr-form}, we obtain
\begin{align*}
  \beta
  &=
  X_0
  \bigl(
    A\omega+\rho_{\mathrm{ss}}(\theta_{\SL_r})
  \bigr)
  X_0^{-1}
  -
  dX_0\,X_0^{-1}\\
  &=
  A\omega
  +
  X_0\rho_{\mathrm{ss}}(\theta_{\SL_r})X_0^{-1}
  -
  X_0\rho_{\mathrm{ss}}(\theta_{\SL_r})X_0^{-1}\\
  &=A\omega.
\end{align*}
Hence, after pullback to the central cover and an algebraic gauge
transformation, the semisimple component disappears and only the
$\Gm$-component remains.

\subsection{Essential surjectivity}

Let $M\in\Invdesc(\GL_r)$. Its pullback $E=p^*M$ belongs to $\Invtriv(\Gm\times\SL_r)$ and carries its canonical $\mu_r$-descent datum. By the preceding gauge argument,
$
  E\cong\pr_{\Gm}^*N
$
for some $N\in\Invtriv(\Gm)$.

Transport the $\mu_r$-linearization through this isomorphism. Under
the deck action
\[
  \zeta\cdot(t,s)=(\zeta t,\zeta^{-1}s),
\]
its structure maps become horizontal isomorphisms
\[
  \pr_{\Gm}^*t_\zeta^*N\longrightarrow\pr_{\Gm}^*N.
\]
Their matrix entries are annihilated by all vector fields tangent to
the connected factor $\SL_r$, and hence depend only on the $\Gm$-coordinate. They are therefore pulled back from isomorphisms
\[
  \beta_\zeta:t_\zeta^*N\isoto N
  \qquad(\zeta\in\mu_r),
\]
which satisfy the descent cocycle condition for the finite \'etale Galois isogeny
\[
  q:\Gm\longrightarrow\Gm,
  \qquad t\longmapsto t^r.
\]
Hence, by finite \'etale descent, $(N,(\beta_\zeta))$ descends to a flat
bundle $\overline N$ on $\Gm$, with a $\mu_r$-equivariant
isomorphism $
  q^*\overline N\cong N$.
By \cref{lem:torus-descent}, there exists
$N_0\in\Invtriv(\Gm)$ such that $\overline N\cong N_0$. Since
\[
  \Det\circ p=q\circ\pr_{\Gm},
\]
pulling back this isomorphism gives a $\mu_r$-equivariant
isomorphism
\[
  p^*M\cong p^*\Det^*N_0.
\]
Finite étale descent therefore yields
\[
  M\cong\Det^*N_0.
\]
We now give the explicit descended gauge needed for the trivial-bundle
category.

Let $M_\rho\in\Invtriv(\GL_r)$ be an invariant connection on the
trivial bundle $\cO_{\GL_r}\otimes V$, corresponding to a Lie
algebra homomorphism $\rho:\mathfrak{gl}_r\longrightarrow\mathfrak{gl}(V)$. Let
\[
  \Phi:\SL_r\longrightarrow\GL(V)
\]
be the algebraic representation integrating
$\rho|_{\mathfrak{sl}_r}$, and put $A=\rho(I_r)$.

The gauge can be made explicit. Since $\mu_r$ is diagonalizable,
the representation
\[
  \mu_r\longrightarrow\GL(V),
  \qquad
  \zeta\longmapsto\Phi(\zeta^{-1}I_r),
\]
gives an isotypic decomposition
\[
  V=\bigoplus_{a\in\mathbb Z/r\mathbb Z}V_a,
  \qquad
  \Phi(\zeta^{-1}I_r)|_{V_a}=\zeta^a\operatorname{id}_{V_a}.
\]
Because $\zeta^{-1}I_r$ is central in $\SL_r$, each $V_a$ is stable under $\Phi(\SL_r)$. Moreover, $A$ commutes with $\Phi(\SL_r)$ and therefore preserves each $V_a$. For every $a$, choose an integer $m_a\equiv a\pmod r$ and define an algebraic homomorphism \[\Psi:\Gm\longrightarrow\GL(V),
  \qquad
  \Psi(t)|_{V_a}=t^{m_a}\operatorname{id}_{V_a}.\]
Then $\Psi(\zeta)=\Phi(\zeta^{-1}I_r)$ for every $\zeta\in\mu_r$. Since $\Psi(t)$ is scalar on each $\Phi(\SL_r)$-stable space $V_a$, it commutes with $\Phi(s)$ for every $s\in\SL_r$; it also commutes with $A$.

Define $X(t,s)=\Phi(s)\Psi(t)^{-1}$. For $\zeta\in\mu_r$, centrality and the defining property of $\Psi$ give 
\begin{align*}
  X(\zeta t,\zeta^{-1}s)
  &=\Phi(\zeta^{-1}I_r)\Phi(s)
    \Psi(\zeta)^{-1}\Psi(t)^{-1}\\
  &=\Phi(s)\Psi(t)^{-1}\\
  &=X(t,s).
\end{align*}

Thus $X$ is $\mu_r$-invariant and descends to an algebraic gauge on $\GL_r$. For completeness, gauge first by $X_0(t,s)=\Phi(s)$, which sends $p^*\alpha$ to $A\omega$. The remaining gauge is $Y=\Psi^{-1}$. Write \[\Psi^{-1}d\Psi=B\omega,
  \qquad B|_{V_a}=m_a\operatorname{id}_{V_a}.\]

Using the gauge convention and the commutation of $A$ with $\Psi$, we obtain \[Y(A\omega)Y^{-1}-(dY)Y^{-1}=(A+B)\omega.\]  
Let $u$ be the coordinate on the target copy of $\Gm$ and set \[N_\rho:=\left(\cO_{\Gm}\otimes V,
  d+\frac{1}{r}(A+B)\frac{du}{u}\right).\] Since $\Det\circ p(t,s)=t^r$, one has $p^*\Det^*\!\left(\frac{du}{u}\right)=r\omega$. Therefore the descended gauge gives $M_\rho\cong\Det^*N_\rho$. Together with the preceding descent argument, this proves essential surjectivity for both categories.

\subsection{Full faithfulness}

We use the following standard quotient calculation.

\begin{lemma}\label{lem:SL-invariants}
For the left action of $\SL_r$ on $\GL_r$, the determinant map is the categorical quotient. Equivalently, pullback induces an isomorphism
\[
  \Det^*:\cO(\Gm)\xrightarrow{\sim}\cO(\GL_r)^{\SL_r}.
\]
\end{lemma}

\begin{proof}
The morphism $\Det:\GL_r\to\Gm$ is a faithfully flat $\SL_r$-torsor. Indeed,
\[
  \SL_r\times\GL_r\isoto\GL_r\times_{\Gm}\GL_r,
  \qquad (h,g)\longmapsto(hg,g).
\]
The equalizer description of faithfully flat descent for functions gives the claim.
\end{proof}

By \cref{prop:MC}, let
$
  M=(\cO_{\Gm}\otimes V,d+A\omega),
  \;
  N=(\cO_{\Gm}\otimes W,d+B\omega),
$
where $\omega=d\log t$. Every object of $\Invdesc(\Gm)$ admits such a presentation by \cref{thm:wei-torus,thm:rs}. A horizontal morphism
$
  f:\Det^*M\longrightarrow\Det^*N
$
is represented by
$
  X\in\cO(\GL_r)\otimes\Hom(V,W)
$
satisfying
\[
  dX+(BX-XA)\,d\log(\Det)=0.
\]
Evaluation on left-invariant vector fields from $\mathfrak{sl}_r$ gives $v(X)=0$, so \cref{lem:SL-invariants} yields a unique
\[
  Y\in\cO(\Gm)\otimes\Hom(V,W)
\]
with $X=\Det^*Y$. Substitution gives
$
  \Det^*\bigl(dY+(BY-YA)\omega\bigr)=0.
$

Apply the algebraic section
\[
  s:\Gm\longrightarrow\GL_r,
  \qquad
  s(t)=\diag(t,1,\dots,1),
\]
which satisfies $\Det\circ s=\operatorname{id}_{\Gm}$. Pulling back
by $s$ gives
$
  dY+(BY-YA)\omega=0.
$
Thus $Y$ is a unique horizontal morphism $M\to N$, proving full faithfulness. The functor is tensor because pullback is tensor, and it is an equivalence by essential surjectivity. This proves the equivalence on $\Invdesc$.

Since the canonical cover of $\Gm$ is the identity,
\[
  \Invtriv(\Gm)=\Invdesc(\Gm).
\]
The essential-surjectivity argument also gives
\[
  \Invtriv(\GL_r)=\Invdesc(\GL_r).
\]
Hence all four categories in \cref{thm:det-equivalence} are abelian by
\cref{thm:rs}, and the two additive equivalences are exact.
This completes the proof of \cref{thm:det-equivalence}.

\section{Connected reductive groups}\label{sec:reductive}

Let $G$ be a connected complex reductive group. Write
\[
  \Gder=[G,G],
  \qquad
  \mathrm{ab}:G\longrightarrow\Gab:=G/\Gder,
\]
and let
$
  \pi_{\mathrm{sc}}:S:=(\Gder)^{\mathrm{sc}}\longrightarrow\Gder,
  \;
  \Gamma:=\ker(\pi_{\mathrm{sc}}).
$
Thus $\Gab$ is a torus and $\Gamma\cong\pi_1((\Gder)^{\mathrm{an}},e)$.

\subsection{Fundamental groups and root data}

\begin{proposition}\label{prop:pi-exact}
There is a natural short exact sequence of finitely generated abelian groups
\[
  0\longrightarrow\Gamma
  \longrightarrow\pi_1(G^{\mathrm{an}},e)
  \xto{\mathrm{ab}_*}
  \pi_1((\Gab)^{\mathrm{an}},e)
  \longrightarrow0.
\]
\end{proposition}

\begin{proof}
The quotient map $G^{\mathrm{an}}\to(\Gab)^{\mathrm{an}}$ is a locally trivial fibration of complex Lie groups with connected fiber $(\Gder)^{\mathrm{an}}$. Since
$(\Gab)^{\mathrm{an}}\cong(\C^*)^d$, we have
$\pi_2((\Gab)^{\mathrm{an}},e)=0$. The homotopy exact sequence therefore gives
\[
  0\to\pi_1((\Gder)^{\mathrm{an}},e)
  \to\pi_1(G^{\mathrm{an}},e)
  \to\pi_1((\Gab)^{\mathrm{an}},e)
  \to0.
\]
Moreover, $S^{\mathrm{an}}\longrightarrow(\Gder)^{\mathrm{an}}$ is the universal covering map, with deck group $\Gamma$, and hence
\[
  \pi_1((\Gder)^{\mathrm{an}},e)\cong\Gamma.
\]
Finally, fundamental groups of topological groups are abelian, and
all the groups above are finitely generated.
\end{proof}

\begin{proposition}[Root-datum description]\label{prop:root-datum}
Let $T\subset G$ be a maximal torus, and let $Q^\vee\subset X_*(T)$ be the lattice generated by the coroots. There is a canonical isomorphism
$
  \pi_1(G^{\mathrm{an}},e)\cong X_*(T)/Q^\vee.
$
\end{proposition}

\begin{proof}
Let $T^{\mathrm{sc}}\subset S=(\Gder)^{\mathrm{sc}}$ be the inverse
image of $T\cap\Gder$. By definition, the algebraic fundamental group
of $G$ is
\[\operatorname{coker}\!\left(
    X_*(T^{\mathrm{sc}})\longrightarrow X_*(T)
  \right).\]
The image of $X_*(T^{\mathrm{sc}})$ in $X_*(T)$ is exactly the
coroot lattice $Q^\vee$. Hence the algebraic fundamental group is
$X_*(T)/Q^\vee$.

The comparison with the topological fundamental group is naturally
written as
$
  \pi_1^{\mathrm{alg}}(G,e)
  \cong
  \pi_1(G^{\mathrm{an}},e)(-1);
$
see \cite[Remark~1.1]{BGA14}. Evaluating the Tate twist (identifying $\mathbb{Z}(1)$ with $\mathbb{Z}$) at the
standard positive generator of $\pi_1(\C^*,e)\cong\mathbb Z$ gives
$
  \pi_1(G^{\mathrm{an}},e)\cong X_*(T)/Q^\vee.
$
\end{proof}
\begin{corollary}\label{cor:root-classification}
For every maximal torus $T\subset G$, there is an exact tensor equivalence
\[
  \Invdesc(G)\isoto\RepC\bigl(X_*(T)/Q^\vee\bigr).
\]
\end{corollary}

\begin{proof}
The result follows from \cref{thm:rs,prop:root-datum}.
\end{proof}

\subsection{Derived monodromy and abelianization}

By \cref{thm:rs}, horizontal sections and monodromy give exact
tensor equivalences
\[
 \operatorname{RH}_{G,e}: \Invdesc(G)\isoto \LocSysC(G^{\mathrm{an}})
  \isoto \RepC\bigl(\pi_1(G^{\mathrm{an}},e)\bigr).
\]

Let $\iota_{\mathrm{der}}:\Gamma\hookrightarrow\pi_1(G^{\mathrm{an}},e)$
be the injection of \cref{prop:pi-exact}. We define the
\emph{derived-monodromy functor} by
\[\operatorname{Res}_{\mathrm{der}}:=\iota_{\mathrm{der}}^*\circ\operatorname{RH}_{G,e}:  \Invdesc(G)\longrightarrow\RepC(\Gamma).
\]
Thus, explicitly,
\[\operatorname{Res}_{\mathrm{der}}(M)=\operatorname{RH}_{G,e}(M)|_\Gamma.
\]
This is an exact tensor functor.

\begin{theorem}[Abelianization and derived monodromy]
\label{thm:abelianization}
Pullback along $\mathrm{ab}:G\to\Gab$ induces a fully faithful exact
tensor functor
\[
  \mathrm{ab}^*:
  \Invdesc(\Gab)\longrightarrow\Invdesc(G).
\]
Its essential image is the full subcategory, consisting precisely of those objects whose
derived monodromy is trivial:
\[
  \left\{
    M\in\Invdesc(G)
    \ \middle|\
    \operatorname{Res}_{\mathrm{der}}(M)
    \cong\mathbf 1_\Gamma^{\oplus\operatorname{rank}(M)}
  \right\}.
\]
Consequently, for every $n\geq1$, pullback induces an injection
\[
  [\Invdesc_n(\Gab)]
  \longrightarrow
  [\Invdesc_n(G)].
\]
These maps are surjective for every $n\geq1$ if and only if
$\Gamma=1$, equivalently, if and only if $\Gder$ is simply
connected.
\end{theorem}
\begin{proof}
Set
\[
  \Pi=\pi_1(G^{\mathrm{an}},e),
  \qquad
  \Lambda=\pi_1((\Gab)^{\mathrm{an}},e).
\]
By \cref{prop:pi-exact}, there is an exact sequence
\[
  0\longrightarrow\Gamma\longrightarrow\Pi
  \longrightarrow\Lambda\longrightarrow0,
\]
so that $\Lambda\cong\Pi/\Gamma$.
By the naturality of monodromy,
under the equivalences of \cref{thm:rs}, pullback along
$\mathrm{ab}$ corresponds to inflation along
$\Pi\twoheadrightarrow\Lambda$:
\[
  \RepC(\Lambda)\longrightarrow\RepC(\Pi),
  \qquad
  \rho\longmapsto\rho\circ\mathrm{ab}_*.
\]
Inflation is a fully faithful exact tensor functor, and its essential
image consists precisely of the representations of $\Pi$ on which
$\Gamma$ acts trivially. This proves the description of the
essential image. Full faithfulness also implies injectivity on
rank-$n$ isomorphism classes for every $n$.

If $\Gamma=1$, then $\Pi\xrightarrow{\sim}\Lambda$, so inflation is
an equivalence. Conversely, suppose that $\Gamma\neq1$. Choose a
nontrivial character
\[
  \chi:\Gamma\longrightarrow\C^*.
\]
Since $\Gab$ is a torus, $\Lambda$ is free abelian. Hence the above
exact sequence of abelian groups in \cref{prop:pi-exact} splits noncanonically, and $\chi$
extends to a character
\[
  \widetilde\chi:\Pi\longrightarrow\C^*.
\]
The corresponding rank-one object of $\Invdesc(G)$ has nontrivial
derived monodromy and therefore does not lie in the essential image
of $\mathrm{ab}^*$. Thus pullback is not surjective even on rank-one
isomorphism classes. This proves the final assertion.
\end{proof}

We have the following corollary. In particular, \cref{conj:wei} holds for $G$.
\begin{corollary}\label{cor:wei-sc-derived}
If $\Gder$ is simply connected, then pullback induces an exact tensor equivalence
\[
  \mathrm{ab}^*: \Invtriv(\Gab)\isoto\Invtriv(G).\]
\end{corollary}

\begin{proof}
Here $\Gamma=1$, so \cref{thm:abelianization} is an equivalence on $\Invdesc$. Since $\Gab$ is a torus, 
\[\Invtriv(\Gab)=\Invdesc(\Gab).\] Pullback of a left-invariant algebraic connection on a trivial vector bundle is again left invariant on a trivial vector bundle. Hence the equivalence restricts to an equivalence on the full subcategories of invariant trivial-bundle connections. Conversely, if $M\in\Invtriv(G)\subset\Invdesc(G)$, essential
surjectivity on $\Invdesc$ gives
\[M\cong\mathrm{ab}^*N\text{ for some }N\in\Invdesc(\Gab)=\Invtriv(\Gab).\qedhere\]
\end{proof}

\subsection{An algebraic construction of derived monodromy}\label{subsec:algebraic-derived-monodromy}
We record an algebraic form of \cref{thm:abelianization}. Put $Z=Z(G)^0$ and define
\[
  \widetilde G=S\times Z,
  \;
  \pi:\widetilde G\longrightarrow G,
  \;
  (s,z)\longmapsto\pi_{\mathrm{sc}}(s)z.
\]
Its kernel
$
  K:=\ker(\pi)
$
is finite and central. The projection $K\to Z$ has image $Z\cap\Gder$ and kernel~$\Gamma$, so
\[K/\Gamma\cong Z\cap\Gder.\]
Moreover, if
$
  q:Z\longrightarrow Z/(Z\cap\Gder)\cong\Gab
$
is the quotient, then
\[
  \mathrm{ab}\circ\pi=q\circ\pr_Z.\]
Let $M\in\Invdesc(G)$. The pullback $\widetilde M=\pi^*M$ is an object of $\Invtriv(\widetilde G)$ with canonical $K$-descent data. Its Lie algebra representation decomposes into commuting representations on $\Lie(S)\oplus\Lie(Z)$. Integrating the semisimple part and applying the corresponding gauge gives an isomorphism 
$ u:\widetilde M \xrightarrow[]{\sim} \pr_Z^*N$ for some $N\in\Invtriv(Z)$.

Transport the $K$-linearization through $u$. Write $k=(s_k,z_k)\in K$ and let \[t_k(s,z)=(ss_k^{-1},zz_k^{-1})\]
be the corresponding right deck transformation. If $\tau_{z_k}:Z\to Z$ is given by $\tau_{z_k}(z)=zz_k^{-1}$, the transported structure map is a horizontal isomorphism \[\lambda_k:t_k^*\pr_Z^*N
  =\pr_Z^*\tau_{z_k}^*N\isoto\pr_Z^*N.\]
Evaluating the horizontal equation in the $\Lie(S)$-directions shows that the matrix entries of $\lambda_k$ are constant along the connected factor $S$. 

Hence there is a unique horizontal isomorphism 
\[
\overline\lambda_k:\tau_{z_k}^*N\isoto N
\qquad\text{with}\qquad
\lambda_k=\pr_Z^*\overline\lambda_k.\]
The maps $\overline\lambda_k$ satisfy the right-action descent
cocycle
\[
  \overline\lambda_{kl}
  =
  \overline\lambda_l\circ
  \tau_{z_l}^*\overline\lambda_k
  \qquad(k,l\in K).
\]
For $\gamma\in\Gamma$, one has $z_\gamma=1$, so $\overline\lambda_\gamma$ is a horizontal automorphism of $N$. Evaluating the right-action cocycle at $1\in Z$ gives $(\overline\lambda_{\gamma_1\gamma_2})_1
  =(\overline\lambda_{\gamma_2})_1
   (\overline\lambda_{\gamma_1})_1$.
Since $\Gamma$ is abelian, interchanging $\gamma_1$ and $\gamma_2$, we therefore obtain a representation, well defined up to conjugacy,
\[
r_M:\Gamma\longrightarrow\GL(N_1),
\qquad r_M(\gamma):=(\overline\lambda_\gamma)_1.
\]
Choose a path $c:[0,1]\to S^{\mathrm{an}}$ from $e$ to $\gamma^{-1}$. Its image under $\pi$ is a loop in $G^{\mathrm{an}}$ representing the image of $\gamma$ under the injection
$\Gamma\hookrightarrow\pi_1(G^{\mathrm{an}},e)$ of \cref{prop:pi-exact}. The lifted path
$t\mapsto(c(t),1)$ has constant $Z$-coordinate, and the connection $\pr_Z^*N$ has no component in the $S$-directions. Hence parallel transport along this lift is the identity on $N_1$. Its endpoint is
$(\gamma^{-1},1)=t_{(\gamma,1)}(e,1)$, so the descent identification of the endpoint fiber with the initial fiber is $(\overline\lambda_\gamma)_1$. Therefore $r_M$ is precisely the derived-monodromy representation $\operatorname{Res}_{\mathrm{der}}(M)$.

This construction is independent of the chosen semisimple gauge. Indeed, if a second gauge gives $u':\widetilde M\isoto\pr_Z^*N'$, then $u'\circ u^{-1}$ is horizontal and, by the same argument in the $S$-directions, is the pullback of a unique horizontal isomorphism $f:N\isoto N'$. The fiber map $f_1$ conjugates the two $\Gamma$-actions. The construction is functorial in horizontal morphisms and compatible with tensor products, and therefore realizes $\operatorname{Res}_{\mathrm{der}}$ algebraically.

Suppose that $\operatorname{Res}_{\mathrm{der}}(M)$ is trivial. Then
$(\overline\lambda_\gamma)_1=\operatorname{id}_{N_1}$ for every
$\gamma\in\Gamma$. A horizontal endomorphism of a connection on the
connected variety $Z$ is determined by its value at one point: after
analytification, it is a morphism of local systems. Hence
$\overline\lambda_\gamma=\operatorname{id}_N$ for every
$\gamma\in\Gamma$. Applying the cocycle with $l=\gamma$, and using
$z_\gamma=1$, gives
\[
  \overline\lambda_{k\gamma}
  =
  \overline\lambda_\gamma\circ
  \tau_{z_\gamma}^*\overline\lambda_k
  =
  \overline\lambda_k.
\]
Thus the transported $K$-linearization depends only on the class of
$k$ in $K/\Gamma$.

Under the identification $K/\Gamma\cong Z\cap\Gder$, these maps are precisely translation descent data on~$N$ for the finite \'etale quotient
\[
  q:Z\longrightarrow Z/(Z\cap\Gder)\cong\Gab.
\]
By finite \'etale descent, these data descend to a flat bundle $\overline N$ on $\Gab$, with an equivariant isomorphism
$
  q^*\overline N\cong N.
$
By \cref{lem:torus-descent}, there exists $N_{\mathrm{ab}}\in\Invtriv(\Gab)$ such that $\overline N\cong N_{\mathrm{ab}}$. Using
$
  \mathrm{ab}\circ\pi=q\circ\pr_Z,
$
we obtain a $K$-equivariant isomorphism
$
  \pi^*M\cong\pi^*\mathrm{ab}^*N_{\mathrm{ab}}.
$
Effective descent along $\pi$ then gives
$
  M\cong\mathrm{ab}^*N_{\mathrm{ab}}.
$
Conversely, the restriction to $\Gamma$ of the descent datum of an object pulled back from $\Gab$ is trivial. This proves algebraically the essential-image statement in \cref{thm:abelianization}.

For completeness, if $\Gamma\neq1$, choose a nontrivial character $\chi:\Gamma\to\C^*$. Since $\Gamma\hookrightarrow K$ is a closed
immersion of finite diagonalizable groups, restriction induces a
surjection
\[
  X^*(K)\twoheadrightarrow X^*(\Gamma).
\]
Thus $\chi$ extends to a character $\widetilde\chi:K\to\C^*$.
Equipping the trivial line connection on $\widetilde G$ with the
$K$-linearization defined by $\widetilde\chi$ and descending along
$\pi$ produces a rank-one object of $\Invdesc(G)$ whose derived monodromy is
$\chi$. In particular, it does not lie in the essential image of
$\mathrm{ab}^*$.

\subsection{The general Karoubi-completion theorem}

\begin{remark}[Monodromy convention]\label{rem:monodromy-convention}
Let $p:Y\to X$ be a connected finite \'etale Galois cover with abelian deck group $\Delta$, and choose basepoints $y_0\in Y$ and $x_0=p(y_0)$. We index the deck transformations by
\[
  t_\delta(y)=y\delta^{-1}\qquad(\delta\in\Delta).
\]
Thus $t_{\delta_1\delta_2}=t_{\delta_1}\circ t_{\delta_2}$. If a linearization is written as maps $A_\delta:t_\delta^*E\to E$, its cocycle identity is
\[
A_{\delta_1\delta_2}
  =A_{\delta_2}\circ t_{\delta_2}^*A_{\delta_1}.
\]
When $\Delta$ is abelian and the $A_\delta$ are constant matrices, this is equivalent to the ordinary representation identity.

The endpoint of a lifted loop defines a quotient homomorphism
\[
  q_p:\pi_1(X,x_0)\twoheadrightarrow
  \pi_1(X,x_0)/p_*\pi_1(Y,y_0)\cong\Delta
\]
by the rule that $q_p([\ell])=\delta$ when the lift of $\ell$ beginning at $y_0$ ends at $t_\delta(y_0)$. 

Suppose that the constant connection with fiber $V$ on $Y$ is
equipped with descent isomorphisms
\[
  A_\delta:
  t_\delta^*(\cO_Y\otimes V,d)
  \isoto
  (\cO_Y\otimes V,d).
\]
Parallel transport along the lifted loop is the identity on $V$,
and the descent identification of the endpoint fiber with the
initial fiber is $(A_\delta)_{y_0}$. Hence the monodromy of
$[\ell]$ is $(A_\delta)_{y_0}$. Thus the descent and monodromy
representations agree with this convention.

If $Y$ is simply connected, then $p$ is the universal cover of $X$,
and hence $q_p$ is an isomorphism: 
\[ q_p:\pi_1(X,x_0)\xrightarrow{\sim}\Delta.
\]
\end{remark}

\begin{theorem}[Karoubi completion for connected reductive groups]\label{thm:general-kar}
For every connected complex reductive group $G$, the natural inclusion induces a tensor equivalence
\[
  \Kar\bigl(\Invtriv(G)\bigr)\isoto\Invdesc(G).
\]
Equivalently, every regular-singular algebraic connection on $G$ is a direct summand of an object of $\Invtriv(G)$.
\end{theorem}

\begin{proof}
By \cref{thm:rs}, it is enough to show that every finite-dimensional representation of
\[
  \Pi:=\pi_1(G^{\mathrm{an}},e)
\]
is a direct summand of the monodromy representation of an object of $\Invtriv(G)$. Let $V$ be such a representation. Since $\Pi$ is abelian and $\Gamma\subset\Pi$ is finite, $V$ decomposes into $\Gamma$-isotypic subrepresentations
$
  V=\bigoplus_{\chi\in\widehat\Gamma}V_\chi,
$
where $\widehat\Gamma:=\Hom(\Gamma,\C^*)$.
Each $V_\chi$ is $\Pi$-stable.

Fix a character $\chi$ occurring in this decomposition. As in the proof of \cref{cor:semisimple-kar}, choose an irreducible algebraic representation
\[
  \Phi_\chi:S=(\Gder)^{\mathrm{sc}}\longrightarrow\GL(W_\chi)\]
whose restriction to $\Gamma$ is the scalar character $\chi$. Let $Z=Z(G)^0$ and let
\[
  \pi:S\times Z\longrightarrow G,
  \;
  K=\ker(\pi),
\]
be the canonical central cover introduced above. Extending $d\Phi_\chi$ by zero on $\Lie(Z)$ and transporting it through the Lie-algebra isomorphism $d\pi$ defines a Lie algebra homomorphism
$
  \Lie(G)\longrightarrow\mathfrak{gl}(W_\chi),
$
and hence an object $E_\chi\in\Invtriv(G)$.

We compute its monodromy. The pullback $\pi^*E_\chi$ is trivialized by the algebraic gauge
$
  X(s,z)=\Phi_\chi(s).
$
Write an element of $K$ as $k=(s_k,z_k)$ and use the deck transformation $t_k(y)=yk^{-1}$ from \cref{rem:monodromy-convention}. Since $K$ is central in $S\times Z$, we have $s_k\in Z(S)$. The canonical $K$-linearization of $\pi^*E_\chi$, transported through the gauge $X$, is the constant matrix
\begin{align*}
  A_k(s,z)
  &=X(s,z)X\bigl(t_k(s,z)\bigr)^{-1}\\
  &=\Phi_\chi(s)\Phi_\chi(ss_k^{-1})^{-1}\\
  &=\Phi_\chi(s_k).
\end{align*}
By Schur's lemma, $\Phi_\chi(s_k)$ is scalar. Thus there is a character
\[\kappa_\chi:K\longrightarrow\C^*,
  \;
\Phi_\chi(s_k)=\kappa_\chi(k)\operatorname{id}_{W_\chi}.\]
After applying the gauge $X$, the connection $\pi^*E_\chi$ is
trivial. Its monodromy is therefore trivial on
\[\pi_*\pi_1((S\times Z)^{\mathrm{an}},e).\]
Hence it factors through the quotient homomorphism
$
  q_\pi:\Pi\twoheadrightarrow K
$
of \cref{rem:monodromy-convention}, and the monodromy of $E_\chi$ is
\[
  \widetilde\chi^{\oplus\dim W_\chi},
  \;
  \widetilde\chi:=\kappa_\chi\circ q_\pi:\Pi\longrightarrow\C^*.
\]
Under the injection $\Gamma\hookrightarrow\Pi$ of
\cref{prop:pi-exact}, an element $\gamma\in\Gamma$ is represented
by a loop whose lift to $S\times Z$ ends at
\[
  (\gamma^{-1},1)=t_{(\gamma,1)}(e,1).
\]
Thus, by \cref{rem:monodromy-convention}, $q_\pi(\gamma)=(\gamma,1)$, and consequently
\[
  \widetilde\chi|_\Gamma=\chi.
\]
Set
\[
  U_\chi:=\widetilde\chi^{-1}\otimes V_\chi.
\]
Since $\widetilde\chi|_\Gamma=\chi$ and $\Gamma$ acts on
$V_\chi$ through $\chi$, the representation $U_\chi$ is trivial
on $\Gamma$. It therefore factors through
\[
  \Pi/\Gamma\cong\pi_1((\Gab)^{\mathrm{an}},e).
\]
By \cref{thm:wei-torus}, this quotient representation is the monodromy of some $N_\chi\in\Invtriv(\Gab)$. Consequently,
\[
  E_\chi\otimes\mathrm{ab}^*N_\chi\in\Invtriv(G)
\]
has monodromy
$
  V_\chi^{\oplus\dim W_\chi}.
$ Indeed, $E_\chi$ contributes the scalar character $\widetilde\chi$ on $W_\chi$, whereas
$\mathrm{ab}^*N_\chi$ contributes the representation $\widetilde\chi^{-1}\otimes V_\chi$.
Thus $V_\chi$ is a direct summand of an object coming from $\Invtriv(G)$. Taking the direct sum over the finitely many characters occurring in $V$ proves the same for $V$.

Thus every object of $\Invdesc(G)$ is a direct summand of an object
of $\Invtriv(G)$. Since $\Invtriv(G)$ is a full tensor subcategory
of the idempotent-complete category $\Invdesc(G)$, the inclusion
extends uniquely to a fully faithful tensor functor
\[
  \Kar\bigl(\Invtriv(G)\bigr)\longrightarrow\Invdesc(G).
\]
The direct-summand statement gives essential surjectivity, and hence
this functor is a tensor equivalence.
\end{proof}

\subsection{Wei's conjecture and regular-singular connections}
{The following result gives a group-theoretic criterion for Wei's conjecture to hold in every rank and, equivalently, for invariant trivial-bundle connections to exhaust all regular-singular connections.}
\begin{theorem}[{The Wei--Regular-Singular Criterion}]\label{cor:wei-iff}
For a connected complex reductive group $G$, the following conditions are equivalent:
\begin{enumerate}
  \item $\Gder$ is simply connected;
  \item pullback along $\mathrm{ab}:G\to\Gab$ induces an equivalence
  \[
    \mathrm{ab}^*:\Invtriv(\Gab)\isoto\Invtriv(G);
  \]
  \item Wei's {\cref{conj:wei}} holds for $G$ in every rank;
  \item $\Invtriv(G)$ is idempotent complete;
  \item $\Invtriv(G)=\Invdesc(G)=\Connrs(G)$.
\end{enumerate}
\end{theorem}

\begin{proof}
Set $\Gamma=\pi_1((\Gder)^{\mathrm{an}},e)$. Then $\Gamma=1$ if and only if $\Gder$ is simply connected.

Suppose first that $\Gamma=1$. By
\cref{cor:wei-sc-derived}, pullback induces an equivalence
\[
  \mathrm{ab}^*:
  \Invtriv(\Gab)\isoto\Invtriv(G),
\]
so (2) holds and hence so does (3). Moreover,
\cref{thm:abelianization} and the torus equality give
\[\Invtriv(G)=\Invdesc(G)=\Connrs(G).
\]
Thus (5) holds. Since $\Connrs(G)$ is idempotent complete, (4)
follows as well. Hence (1) implies all the remaining conditions.

Conversely, suppose that $\Gamma\neq1$, and choose a nontrivial
character $\chi:\Gamma\longrightarrow\C^*$. The construction in the proof of \cref{thm:general-kar} produces
an object $E_\chi\in\Invtriv(G)$ whose derived monodromy is a
nonzero direct sum of copies of $\chi$. By
\cref{thm:abelianization}, $E_\chi$ does not lie in the essential
image of $\mathrm{ab}^*$. Thus (2) fails, and the same object gives
a counterexample to Wei's \cref{conj:wei}; hence (3) fails.

We next show that (4) fails. The construction in the proof of
\cref{thm:abelianization} gives a rank-one object
\[
  L_\chi\in\Invdesc(G)
\]
with derived monodromy $\chi$. On the other hand, every rank-one
object of $\Invtriv(G)$ has trivial derived monodromy. Indeed, such
an object is determined by a Lie algebra character
\[
  \rho:\Lie(G)\longrightarrow\C,
\]
which vanishes on $[\Lie(G),\Lie(G)]=\Lie(\Gder)$.
After pullback to the canonical cover, its connection is therefore
pulled back from the central torus in the canonical trivialization,
so the restricted $\Gamma$-linearization is trivial. Consequently,
\[
  L_\chi\notin\Invtriv(G).
\]
By \cref{thm:general-kar}, however, $L_\chi$ is a direct summand in
$\Invdesc(G)$ of an object of $\Invtriv(G)$. Since
$\Invtriv(G)$ is a full subcategory of $\Invdesc(G)$, the
corresponding idempotent is a morphism in $\Invtriv(G)$ but does
not split there. Thus $\Invtriv(G)$ is not idempotent complete, so~(4) fails. In particular, (5) also fails.

Therefore, when $\Gamma\neq1$, none of (1)--(5) holds. This proves
their equivalence.
\end{proof}

\begin{corollary}\label{cor:simples-semisimple}
Every simple object of $\Invdesc(G)$ has rank one. Moreover, this category is semisimple if and only if $G$ is semisimple.
\end{corollary}

\begin{proof}
Set $\Pi=\pi_1(G^{\mathrm{an}},e)$. By \cref{thm:rs}, there is an exact equivalence
\[
  \Invdesc(G)\simeq\RepC(\Pi),
\]
under which rank corresponds to vector-space dimension. Since
$\Pi$ is abelian, every finite-dimensional irreducible complex
representation of $\Pi$ is one-dimensional. Hence every simple
object of $\Invdesc(G)$ has rank one.

If $G$ is semisimple, then $\Pi$ is finite, so Maschke's theorem
implies that $\RepC(\Pi)$, and therefore $\Invdesc(G)$, is
semisimple.

Conversely, suppose that $G$ is not semisimple. Then $\Gab$ has
positive dimension, and hence
\[
  \pi_1((\Gab)^{\mathrm{an}},e)\cong\mathbb Z^d\qquad (d>0).
\]
By \cref{prop:pi-exact}, $\Pi$ surjects onto $\mathbb Z^d$, and
therefore onto $\mathbb Z$. Consider the two-dimensional
representation of $\mathbb Z$ in which a generator acts by
\[
  \begin{pmatrix}
    1&1\\
    0&1
  \end{pmatrix}.
\]
It is a nonsplit extension of the trivial representation by itself.
Inflating it along $\Pi\twoheadrightarrow\mathbb Z$ gives a
nonsplit extension in $\RepC(\Pi)$. Thus neither $\RepC(\Pi)$ nor
$\Invdesc(G)$ is semisimple.
\end{proof}

\subsection{Example: \texorpdfstring{$\Gm\times\PGL_2$}{Gm x PGL2}}

Let
$
  G=\Gm\times\PGL_2.
$
Then
$
  \pi_1(G^{\mathrm{an}},e)\cong\mathbb Z\times\mu_2,
$
so
\[
  \Invdesc(G)\simeq\RepC(\mathbb Z\times\mu_2).
\]
The essential image of pullback from the abelianization $\Gm$ is the full subcategory on which $\mu_2$ acts trivially.

For a concrete family, let $L\in\Invtriv_1(\Gm)$ and let $V$ be a finite-dimensional representation of~$\mu_2$. On
$
  \widetilde G=\Gm\times\SL_2
$
consider
\[
  \pr_{\Gm}^*L\otimes(\cO_{\widetilde G}\otimes V,d),
\]
with $\mu_2$ acting by deck transformations on $\SL_2$ and through its given action on $V$. It descends to an object $M_{L,V}\in\Invdesc(G)$ whose derived monodromy is $V$. It lies in the essential image of pullback from the abelianization precisely when the $\mu_2$-action on $V$ is trivial.

\section{De Rham cohomology in the semisimple case}\label{sec:derham}

Throughout this section, $G$ is connected semisimple,
$\pi:\Gsc\to G$ is the simply connected central cover, and
$\Gamma=\ker(\pi)$. Let $V\in\RepC(\Gamma)$, and write
$\rho_V:\Gamma\to\GL(V)$ for the corresponding representation.

Let 
$
  \widetilde M_V:=(\cO_{\Gsc}\otimes_{\C}V,d).
$ For $\gamma\in\Gamma$, let
$t_\gamma(s)=s\gamma^{-1}$ be the right deck transformation fixed in
\cref{rem:monodromy-convention}.  We equip $\widetilde M_V$ with the
$\Gamma$-linearization whose descent map
\[
  A_\gamma:t_\gamma^*\widetilde M_V\longrightarrow\widetilde M_V
\]
is the constant horizontal automorphism $\rho_V(\gamma)$ on the
factor $V$.  The descent cocycle is $A_{\gamma_1\gamma_2}
  =A_{\gamma_2}\circ t_{\gamma_2}^*A_{\gamma_1}$. Because the $A_\gamma$ are constant, $\Gamma$ is abelian, and $\rho_V$ is a representation, this identity is satisfied. Finite \'{e}tale descent therefore produces an algebraic
connection
$
  M_V:=\bigl(\pi_*\widetilde M_V\bigr)^\Gamma
  \;\text{on }G,
$
and the canonical pullback isomorphism
$\pi^*M_V\cong\widetilde M_V$ is $\Gamma$-equivariant.  In particular,
$M_V\in\Invdesc(G)$.

\begin{remark}\label{rem:MV-two-categories}
For every $V\in\RepC(\Gamma)$, the construction above gives $M_V\in\Invdesc(G)$. By \cref{prop:semisimple-triv}, 
$M_V\in\Invtriv(G)$ precisely when the isomorphism class of
$\rho_V$ lies in the image of
\[
  \frac{\Hom_{\mathrm{alg.grp}}(\Gsc,\GL(V))}{\GL(V)}
  \longrightarrow
  \frac{\Hom(\Gamma,\GL(V))}{\GL(V)}.
\]   
\end{remark}

For an algebraic connection $E=(\mathcal E,\nabla)$, write
\[
  \DR(E):=
  \bigl[
    \mathcal E\xto{\nabla}
    \mathcal E\otimes\Omega^1_{G/\C}\longrightarrow
    \mathcal E\otimes\Omega^2_{G/\C}\longrightarrow\cdots
  \bigr]
\]
for its algebraic de Rham complex, and put
\[
  \RGammadR(G,E):=R\Gamma(G,\DR(E)).
\]

\begin{theorem}[De Rham cohomology of descended connections]\label{thm:derham}
There is a canonical isomorphism in the derived category of complex vector spaces
\[
  \RGammadR(G,M_V)
  \cong
  \bigl(
    \RGammadR(\Gsc,\cO_{\Gsc})\otimes_{\C}V
  \bigr)^\Gamma.
\]
Consequently, for every $i\geq0$,
\[
  \HdR^i(G,M_V)
  \cong
  \bigl(\HdR^i(\Gsc,\C)\otimes_{\C}V\bigr)^\Gamma
  \cong
  \HdR^i(\Gsc,\C)\otimes_{\C}V^\Gamma.
\]
In particular,
$
  \HdR^0(G,M_V)\cong V^\Gamma.
$
\end{theorem}

\begin{proof}
\noindent\emph{Step 1: the de Rham complex descends term by term.}
Write $\widetilde E_V=\cO_{\Gsc}\otimes_{\C}V$ for the underlying
bundle of $\widetilde M_V$.  Effective descent gives
$
  E_V\cong(\pi_*\widetilde E_V)^\Gamma,
$
where $E_V$ is the underlying bundle of $M_V$.  Because $\pi$ is
\'{e}tale, pullback of differentials gives canonical isomorphisms
\[
  \Omega^j_{\Gsc/\C}\cong\pi^*\Omega^j_{G/\C}
  \qquad(j\geq0).
\]
The projection formula therefore yields
\[
  \pi_*\bigl(\widetilde E_V\otimes\Omega^j_{\Gsc/\C}\bigr)
  \cong
  \pi_*\widetilde E_V\otimes\Omega^j_{G/\C}.
\]
These isomorphisms are compatible with the de Rham differentials and
yield a $\Gamma$-equivariant isomorphism of complexes
\[
  \pi_*\DR(\widetilde M_V)
  \cong
  (\pi_*\widetilde E_V)
  \otimes_{\cO_G}\Omega^\bullet_{G/\C}.
\]
Since $\Gamma$ is finite and the ground field has characteristic zero,
the averaging operator
\[
  e_\Gamma=\frac{1}{|\Gamma|}\sum_{\gamma\in\Gamma}\gamma
\]
is an idempotent whose image is the invariant subsheaf.  Thus taking
$\Gamma$-invariants is exact and commutes with tensoring by the locally
free sheaves $\Omega^j_{G/\C}$.  Consequently,
\begin{equation}\label{eq:DR-descent-det}
  \DR(M_V)
  \cong
  \bigl(\pi_*\DR(\widetilde M_V)\bigr)^\Gamma.
\end{equation}
\smallskip
\noindent\emph{Step 2: passage to derived global sections.}
Because the invariant complex is the image of the idempotent
$e_\Gamma$, applying the additive functor $R\Gamma(G,-)$ preserves this
splitting.  Hence \eqref{eq:DR-descent-det} gives
\begin{align*}
  \RGammadR(G,M_V)
  &\cong
  R\Gamma\bigl(G,(\pi_*\DR(\widetilde M_V))^\Gamma\bigr)\\
  &\cong
  \bigl(R\Gamma(G,\pi_*\DR(\widetilde M_V))\bigr)^\Gamma.
\end{align*}
A finite morphism is affine.  For a complex of quasi-coherent sheaves
on the source of an affine morphism, ordinary pushforward computes
derived pushforward. Therefore
$
  R\pi_*\DR(\widetilde M_V)
  \cong
  \pi_*\DR(\widetilde M_V).
$
The Leray identity for the composition with the structure morphism
then gives
\[
  R\Gamma(G,\pi_*\DR(\widetilde M_V))
  \cong
  R\Gamma(\Gsc,\DR(\widetilde M_V)).
\]
Finally, $\widetilde M_V$ is the trivial connection with constant fiber
$V$, so there is a canonical identification of complexes
$
  \DR(\widetilde M_V)
  \cong
  \DR(\cO_{\Gsc})\otimes_{\C}V.
$
Combining these isomorphisms proves
\[
  \RGammadR(G,M_V)
  \cong
  \bigl(
    \RGammadR(\Gsc,\cO_{\Gsc})\otimes_{\C}V
  \bigr)^\Gamma.
\]

\smallskip
\noindent\emph{Step 3: passage to cohomology.}
Exactness of invariants implies
$H^i(C^\Gamma)\cong H^i(C)^\Gamma$ for every
$\Gamma$-equivariant complex $C$.  Since tensoring over $\C$ with the
finite-dimensional vector space $V$ is exact, the derived formula
therefore gives
\[
  \HdR^i(G,M_V)
  \cong
  \bigl(\HdR^i(\Gsc,\C)\otimes_{\C}V\bigr)^\Gamma.
\]

\smallskip
\noindent\emph{Step 4: the deck group acts trivially on
$\HdR^*(\Gsc,\C)$.}
Let
$m:\Gsc\times\Gsc\to\Gsc$ be multiplication, let
$e:\Spec(\C)\to\Gsc$ be the identity, and let
$g:\Spec(\C)\to\Gsc$ be any complex point. The de Rham K\"unneth
isomorphism identifies
\[
  \HdR^*(\Gsc\times\Gsc,\C)
  \cong
  \HdR^*(\Gsc,\C)\otimes_{\C}\HdR^*(\Gsc,\C).
\]
Since $L_g=m\circ(g\times\operatorname{id})$, functoriality gives
$
  L_g^*=(g^*\otimes\operatorname{id})\circ m^*.
$
For positive degrees, both $g^*$ and $e^*$ take values in
$\HdR^{>0}(\Spec(\C),\C)=0$.  In degree zero they agree because
$\Gsc$ is connected and $\HdR^0(\Gsc,\C)=\C$.  Hence $g^*=e^*$ in
all degrees, and the counit identity 
gives
\[
  L_g^*
  =(e^*\otimes\operatorname{id})\circ m^*
  =\operatorname{id}.
\]
The deck transformation $t_\gamma$ is right translation by
$\gamma^{-1}$; since $\Gamma\subset Z(\Gsc)$, it is also left
translation by $\gamma^{-1}$.  Thus every $\gamma\in\Gamma$ acts
trivially on $\HdR^i(\Gsc,\C)$.  It follows that
\[
  \bigl(\HdR^i(\Gsc,\C)\otimes_{\C}V\bigr)^\Gamma
  =
  \HdR^i(\Gsc,\C)\otimes_{\C}V^\Gamma.
\]
Taking $i=0$ proves $\HdR^0(G,M_V)\cong V^\Gamma$.
\end{proof}

\begin{corollary}\label{cor:cohomology-cover}
Pullback along $\pi$ induces canonical isomorphisms
\[
  \RGammadR(G,\cO_G)\cong\RGammadR(\Gsc,\cO_{\Gsc}),
  \qquad
  \HdR^i(G,\C)\cong\HdR^i(\Gsc,\C).
\]
\end{corollary}

\begin{proof}
Take $V=\C$ with the trivial $\Gamma$-action.  Its descent is the
trivial connection $(\cO_G,d)$, and \cref{thm:derham} identifies
$
  \RGammadR(G,\cO_G)
  \cong
  \RGammadR(\Gsc,\cO_{\Gsc})^\Gamma.
$
Under this identification, pullback along~$\pi$ is the canonical map
from the invariant direct summand to the full complex.  On cohomology
this map is
$
  \HdR^i(\Gsc,\C)^\Gamma
  \longrightarrow
  \HdR^i(\Gsc,\C).
$
Step~4 in the proof of \cref{thm:derham} shows that the $\Gamma$-action
on the target is trivial, so this map is an isomorphism for every $i$.
Hence the map of derived de Rham complexes is a quasi-isomorphism,
which proves both assertions.
\end{proof}

\begin{corollary}[Borel's description]\label{cor:borel}
Let $d_1,\dots,d_r$ be the fundamental degrees of the Weyl group of $G$, where $r=\operatorname{rank}(G)$. There is a noncanonical isomorphism of graded algebras
\[
  \HdR^*(G,\C)
  \cong
  \Lambda_{\C}(x_1,\dots,x_r),
  \qquad
  \deg(x_j)=2d_j-1.
\]
For $V\in\RepC(\Gamma)$,
\[
  \HdR^*(G,M_V)
  \cong
  \Lambda_{\C}(x_1,\dots,x_r)\otimes_{\C}V^\Gamma
\]
as graded vector spaces.
\end{corollary}

\begin{proof}
Let $K\subset(\Gsc)^{\mathrm{an}}$ be a maximal compact subgroup.
By polar decomposition, $K$ is a deformation retract of
$(\Gsc)^{\mathrm{an}}$. Borel's theorem
\cite[\S27, Proposition~27.2]{Bor53} identifies
$H^*(K,\C)$ with an exterior algebra on generators of degrees
$2d_j-1$.  The algebraic de Rham comparison isomorphism gives
\[
  \HdR^*(\Gsc,\C)
  \cong
  H^*((\Gsc)^{\mathrm{an}},\C)
  \cong
  H^*(K,\C).
\]
Combining this with \cref{cor:cohomology-cover} proves the first
statement.  The second is the cohomology formula in
\cref{thm:derham} together with the triviality of the $\Gamma$-action
on $\HdR^*(\Gsc,\C)$.
\end{proof}

\section{Local systems and monodromy in the semisimple case}\label{sec:local-systems}

We continue with the notation of \cref{sec:derham} and write
$M_V=(E_V,\nabla_V)$. The analytification $E_V^{\mathrm{an}}$ is a holomorphic vector bundle
on the complex manifold $G^{\mathrm{an}}$, and $\nabla_V$ induces an
integrable holomorphic connection
\[
  \nabla_V^{\mathrm{an}}:
  E_V^{\mathrm{an}}
  \longrightarrow
  E_V^{\mathrm{an}}\otimes_{\cO_{G^{\mathrm{an}}}}
  \Omega^1_{G^{\mathrm{an}}}.
\]
We define
\[
  \mathcal L_V:=\ker(\nabla_V^{\mathrm{an}})
\]
in the analytic topology.  Flatness implies that $\mathcal L_V$ is a
complex local system of rank $\dim_{\C}V$.

\begin{lemma}\label{lem:universal-cover}
The analytification
$
  \pi^{\mathrm{an}}:(\Gsc)^{\mathrm{an}}\longrightarrow G^{\mathrm{an}}
$
is the universal covering map, with deck group $\Gamma$.
\end{lemma}

\begin{proof}
The central isogeny $\pi$ is finite and has finite reduced kernel
$\Gamma$; over $\C$ it is therefore finite \'{e}tale.  By the Riemann
existence theorem, analytification sends a finite \'{e}tale morphism to
a finite topological covering; see
\cite[Expos\'{e}~XII, Th\'{e}or\`eme~5.1]{GR71}.  Thus
$\pi^{\mathrm{an}}$ is a covering map of connected complex manifolds.

The right action
\[
  t_\gamma(s)=s\gamma^{-1}
  \qquad(\gamma\in\Gamma)
\]
is free, preserves every fiber of $\pi^{\mathrm{an}}$, and acts simply
transitively on each fiber.  Hence these transformations are exactly
the deck transformations, and the deck group is $\Gamma$.
Finally, $(\Gsc)^{\mathrm{an}}$ is simply connected
\cite[Appendix~D, Proposition~D.4.1]{Con14}.  A connected covering
with simply connected total space is universal, proving the lemma.
\end{proof}

\begin{proposition}[Monodromy]\label{prop:monodromy}
With the convention of \cref{rem:monodromy-convention}, the monodromy of $\mathcal L_V$ is the composite
\[
  \pi_1(G^{\mathrm{an}},e)\isoto\Gamma
  \longrightarrow\GL(V),
\]
where the second map is the given representation. In particular, every finite-rank local system on $G^{\mathrm{an}}$ has finite monodromy.
\end{proposition}

\begin{proof}
By construction of $M_V$, there is a $\Gamma$-equivariant
isomorphism
\[
  \pi^*M_V\cong(\cO_{\Gsc}\otimes_{\C}V,d).
\]
After analytification and passage to horizontal sections, this says
that $(\pi^{\mathrm{an}})^{-1}\mathcal L_V$ is the constant local
system $\underline V$ on $(\Gsc)^{\mathrm{an}}$.  Its descent map for
$\gamma\in\Gamma$ is the constant linear transformation
$\rho_V(\gamma)$.

Let $\ell$ be a loop in $G^{\mathrm{an}}$ based at $e$, and lift it to
a path $\widetilde\ell$ beginning at the identity
$\widetilde e\in(\Gsc)^{\mathrm{an}}$.  Since the covering is
universal, there is a unique $\gamma\in\Gamma$ such that
\[
  \widetilde\ell(1)=t_\gamma(\widetilde e).
\]
By \cref{rem:monodromy-convention}, the assignment
$[\ell]\mapsto\gamma$ is the isomorphism
$\pi_1(G^{\mathrm{an}},e)\isoto\Gamma$.  Parallel transport in the
constant local system $\underline V$ carries $v\in V$ to the same
vector $v$ at the endpoint of $\widetilde\ell$.  To identify that
endpoint fiber with the fiber over $\widetilde e$, one applies the
descent map $A_\gamma=\rho_V(\gamma)$.  Therefore the monodromy of
$[\ell]$ is exactly $\rho_V(\gamma)$.

Finally, every finite-rank local system on the connected manifold
$G^{\mathrm{an}}$ is determined by a finite-dimensional
representation of its fundamental group.  Since
$\pi_1(G^{\mathrm{an}},e)\cong\Gamma$ is finite, the image of every
such representation is finite.
\end{proof}

\begin{theorem}[Classification of local systems]\label{thm:local-systems}
The assignment
\[
  V\longmapsto\mathcal L_V
\]
defines an exact tensor equivalence
\[
  \RepC(\Gamma)\isoto\LocSysC(G^{\mathrm{an}}).
\]
Thus every local system has finite monodromy, the category is semisimple, and its simple objects are rank-one local systems indexed by the characters of $\Gamma$.
\end{theorem}

\begin{proof}
For a connected, locally path-connected, and semilocally simply
connected space $X$ with basepoint $x$, the monodromy functor
\[
  \operatorname{Mon}_x:\LocSysC(X)
  \longrightarrow
  \RepC(\pi_1(X,x))
\]
sends a local system to its fiber at $x$ equipped with parallel
transport around loops.  It is an equivalence: a quasi-inverse sends
a representation $\rho$ on $W$ to the local system obtained by
descending the constant local system $\underline W$ on the universal
cover, with the deck group acting on $W$ through $\rho$.

Apply this to $X=G^{\mathrm{an}}$.  By \cref{lem:universal-cover}, its
universal cover is $(\Gsc)^{\mathrm{an}}$ and its deck group is
$\Gamma$.  By \cref{prop:monodromy}, the quasi-inverse just described
sends $V\in\RepC(\Gamma)$ to the specific descended local system
$\mathcal L_V$.  Hence $V\mapsto\mathcal L_V$ is fully faithful and
essentially surjective.  Exactness follows because the functor is
computed on the fiber at the basepoint, and tensor compatibility
follows from
\[
  \operatorname{Mon}_x(\mathcal L\otimes\mathcal K)
  =\operatorname{Mon}_x(\mathcal L)\otimes
   \operatorname{Mon}_x(\mathcal K).
\]
Since $\Gamma$ is finite, Maschke's theorem makes
$\RepC(\Gamma)$ semisimple.  Moreover, $\Gamma$ is abelian because it
is contained in the center of $\Gsc$, so every irreducible complex
representation is one-dimensional and is given by a character of
$\Gamma$.
\end{proof}

\begin{theorem}[Betti cohomology]\label{thm:betti}
For every $i\geq0$, there is a canonical isomorphism
\[
  H^i(G^{\mathrm{an}},\mathcal L_V)
  \cong
  \bigl(H^i((\Gsc)^{\mathrm{an}},\C)\otimes_{\C}V\bigr)^\Gamma.
\]
Under the regular-singular de Rham--Betti comparison, this identification agrees with the formula in \cref{thm:derham}:
$
  \HdR^i(G,M_V)\isoto H^i(G^{\mathrm{an}},\mathcal L_V).
$
\end{theorem}

\begin{proof}
Set $p=\pi^{\mathrm{an}}$ and let $\underline V$ denote the constant
sheaf with fiber $V$ on $(\Gsc)^{\mathrm{an}}$, equipped with the
$\Gamma$-linearization used above.  Descent gives
\begin{equation}\label{eq:LV-analytic-descent}
  \mathcal L_V\cong(p_*\underline V)^\Gamma.
\end{equation}
Indeed, on an evenly covered open set $U\subset G^{\mathrm{an}}$, the
inverse image $p^{-1}(U)$ is the finite disjoint union of copies of
$U$, and the invariant sections of the corresponding product of
constant sheaves are precisely the sections of the descended local
system.

Taking invariants is exact because the averaging idempotent splits the inclusion of the invariant subspace on every stalk. The direct-image functor $p_*$ is also exact for a
finite covering: on an evenly covered open set it is a finite product
of restriction functors, and finite products of complex vector spaces
are exact.  Moreover,
\[
  \Gamma(G^{\mathrm{an}},p_*\mathcal F)
  =
  \Gamma((\Gsc)^{\mathrm{an}},\mathcal F)
\]
for every sheaf $\mathcal F$ on $(\Gsc)^{\mathrm{an}}$.  Equivalently,
$p_*$ preserves injectives because it is right adjoint to the exact
inverse-image functor.  Therefore \eqref{eq:LV-analytic-descent}
gives canonical derived isomorphisms
\begin{align*}
  R\Gamma(G^{\mathrm{an}},\mathcal L_V)
  &\cong
  \bigl(R\Gamma(G^{\mathrm{an}},p_*\underline V)\bigr)^\Gamma\\
  &\cong
  \bigl(R\Gamma((\Gsc)^{\mathrm{an}},\underline V)\bigr)^\Gamma\\
  &\cong
  \bigl(R\Gamma((\Gsc)^{\mathrm{an}},\C)
         \otimes_{\C}V\bigr)^\Gamma.
\end{align*}
The last isomorphism uses that $V$ is finite-dimensional, so tensoring
with $V$ is exact.  Taking degree-$i$ cohomology proves
\[
  H^i(G^{\mathrm{an}},\mathcal L_V)
  \cong
  \bigl(H^i((\Gsc)^{\mathrm{an}},\C)\otimes_{\C}V\bigr)^\Gamma.
\]
Each deck transformation is a translation on the connected complex
Lie group $(\Gsc)^{\mathrm{an}}$ and is homotopic to the identity.
Thus $\Gamma$ acts trivially on
$H^i((\Gsc)^{\mathrm{an}},\C)$. Hence the right-hand side is
$H^i((\Gsc)^{\mathrm{an}},\C)\otimes_{\C}V^\Gamma$.

By \cref{rem:MV-two-categories}, we have  $M_V\in\Invdesc(G)$, and \cref{thm:rs} identifies this category with
$\Connrs(G)$.  Hence $M_V$ is regular-singular. By Deligne's comparison theorem \cite[II.6.2]{Del70}, we have the functorial isomorphism
\[
  \HdR^i(G,M_V)
  \isoto
  H^i(G^{\mathrm{an}},\mathcal L_V).\qedhere
\]
\end{proof}

\end{document}